\documentclass{article}

\usepackage{templates/arxiv}
\usepackage{natbib}
\setcitestyle{authoryear,round,citesep={;},aysep={,},yysep={;}}

\usepackage{amsmath,amsfonts,bm}

\def\eqref#1{equation~\ref{#1}}

\def\1{\bm{1}}

\def\rvx{{\mathbf{x}}}

\def\vzero{{\bm{0}}}
\def\vone{{\bm{1}}}

\def\vtheta{{\bm{\theta}}}
\def\vlambda{{\bm{\lambda}}}

\def\vb{{\bm{b}}}
\def\vc{{\bm{c}}}
\def\vd{{\bm{d}}}

\def\vg{{\bm{g}}}
\def\vh{{\bm{h}}}

\def\vr{{\bm{r}}}
\def\vs{{\bm{s}}}

\def\vu{{\bm{u}}}
\def\vv{{\bm{v}}}
\def\vw{{\bm{w}}}
\def\vx{{\bm{x}}}
\def\vy{{\bm{y}}}
\def\vz{{\bm{z}}}

\def\mA{{\bm{A}}}
\def\mB{{\bm{B}}}
\def\mC{{\bm{C}}}

\def\mE{{\bm{E}}}

\def\mG{{\bm{G}}}
\def\mH{{\bm{H}}}
\def\mI{{\bm{I}}}

\def\mR{{\bm{R}}}

\def\mU{{\bm{U}}}
\def\mV{{\bm{V}}}

\DeclareMathAlphabet{\mathsfit}{\encodingdefault}{\sfdefault}{m}{sl}
\SetMathAlphabet{\mathsfit}{bold}{\encodingdefault}{\sfdefault}{bx}{n}

\newcommand{\R}{\mathbb{R}}

\DeclareMathOperator*{\argmin}{arg\,min}

\usepackage{mathtools}
\usepackage{amssymb}

\renewcommand{\eqref}[1]{\textup{(\ref{#1})}}

\DeclarePairedDelimiter{\abs}{\lvert}{\rvert}
\DeclarePairedDelimiter{\norm}{\lVert}{\rVert}

\newcommand{\ip}[2]{\left\langle #1,\,#2\right\rangle}

\newcommand{\oracle}{\mathsf{O}}

\DeclareMathOperator{\dist}{dist}
\DeclareMathOperator{\dom}{dom}

\DeclareMathOperator{\Lip}{Lip}
\DeclareMathOperator{\prox}{prox}

\def\vxi{{\bm{\xi}}}
\def\vnu{{\bm{\nu}}}

\usepackage{amsthm}

\newtheorem{theorem}{Theorem}[section]
\newtheorem{lemma}[theorem]{Lemma}
\newtheorem{corollary}[theorem]{Corollary}
\newtheorem{proposition}[theorem]{Proposition}
\theoremstyle{definition}
\newtheorem{definition}[theorem]{Definition}

\usepackage{float}
\usepackage[ruled,vlined,linesnumbered]{algorithm2e}
\usepackage{placeins}
\SetAlgoNlRelativeSize{0}
\SetAlgoRefRelativeSize{0}
\SetKwInput{KwIn}{Input}
\SetKwInput{KwInit}{Initialize}
\SetKwInput{KwOut}{Output}
\DontPrintSemicolon

\usepackage{longtable}
\usepackage{array}

\usepackage{hyperref}
\usepackage{url}
\usepackage{booktabs}

\title{Optimal First-Order Methods for a Class of Structured H\"older-Smooth Problems}

\newcommand{\CtrlTotalInstances}{250}

\newcommand{\CtrlLowPHighKCount}{60}
\newcommand{\CtrlVvsUFGMRatio}{21.8}
\newcommand{\CtrlVvsGAGDRatio}{103.1}
\newcommand{\CtrlVvsAGDARatio}{134.5}
\newcommand{\CtrlVModerateHits}{201}

\newcommand{\CtrlUFGMStrictHits}{128}

\newcommand{\RDOLSLooseHits}{17}

\newcommand{\RDBTStrictHits}{15}
\newcommand{\RDGAGDStrictHits}{13}
\newcommand{\RDUFGMStrictHits}{14}

\newcommand{\RDOLSStrictHits}{16}
\newcommand{\RDProblemCount}{18}


\author{%
  Baiying Xin\\ \email{by_xin@sjtu.edu.cn}
  \And
  Qi Deng\\ \email{qdeng24@sjtu.edu.cn}%
}
\affiliation{%
  Antai College of Economics and Management\\
  Shanghai Jiao Tong University, Shanghai, China%
}
\date{}

\begin{document}
\maketitle

\begin{abstract}
    Convex optimization with H\"older-continuous gradients has attracted considerable attention as a framework bridging smooth and nonsmooth optimization. For this class, Nesterov's universal fast gradient method attains the optimal black-box complexity, and we study whether exploiting further problem structure can further improve this bound. In this paper, we study an important class of H\"older-smooth problems, namely regularized $\ell_p$ regression with $p\in(1,2)$, whose loss has a $(p-1)$-H\"older-continuous gradient and for which the universal method requires $\mathcal{O}(\varepsilon^{-2/(3p-2)})$ iterations. Exploiting the problem structure, we characterize a family of smooth approximations based on $p$-Huber smoothing and prove that, for any prescribed smoothness level, it minimizes the worst-case approximation error among smooth lower approximations of the scalar loss. Combined with accelerated proximal gradients, the smoothing approach attains an $\varepsilon$-optimal solution in $\mathcal{O}(\varepsilon^{-1/p})$ matrix--vector products, which significantly improves the complexity guarantee for $1<p<2$. We show that this bound is optimal by providing a matching lower bound. Numerical experiments on synthetic and real data demonstrate the practical performance of the proposed methods.
\end{abstract}

\section{Introduction}

Consider the following composite \(\ell_p\) regression problem for \(p\in[1,2)\):
\begin{equation}\label{eq:lp-regression}
    \min_{\vu\in\R^n} \left\{
        F(\vu) := f(\vu) + \omega(\vu) =
        \frac{1}{p}\norm{\mA\vu-\vb}_p^p + \omega(\vu)
    \right\},
    \quad \mA\in\R^{m\times n}, \vb\in\R^m.
\end{equation}
Here, \(\omega\) represents additional regularization or constraints, such as sparsity-promoting \(\ell_1\) penalties or box constraints. For \(p\in(1,2)\), the loss interpolates between least absolute deviations and least squares. Such \(\ell_p\) regression has long been studied as an alternative to least squares under non-Gaussian errors \citep{money1982linear}. Related objectives also arise in semi-supervised learning and local graph clustering \citep{liu2020strongly,el2016asymptotic}.

For general H\"older-smooth problems, the universal fast gradient method of \citet{nesterov2015universal} achieves an \(\varepsilon\)-optimal solution in \(\mathcal O(\varepsilon^{-2/(3p-2)})\) iterations, without prior knowledge of the H\"older exponent or constant. This rate is optimal in the standard first-order oracle model, so any improvement must exploit additional problem structure. A classical example is \(\ell_1\) regression: while the general nonsmooth bound is \(\mathcal O(\varepsilon^{-2})\), smoothing the absolute-value loss improves the iteration complexity to \(\mathcal O(\varepsilon^{-1})\) \citep{nesterov2005smooth,beck2012smoothing}. This naturally raises the question of whether a similar improvement is possible for \(\ell_p\) regression when \(1<p<2\).

A natural starting point is the Moreau envelope of the scalar loss \(r \mapsto\abs{r}^p/p\), obtained through a quadratic proximal penalty \citep{beck2012smoothing}. At \(p=1\), this gives the classical Huber approximation with a uniformly bounded error. For \(1<p<2\), however, the gap between the scalar loss and its Moreau envelope is unbounded as the residual magnitude increases, for any fixed positive smoothing parameter (Proposition~\ref{prop:moreau-gap}). We adopt the smoothing construction of \citet{bubeck2018homotopy}, which we call \(p\)-Huber smoothing (Definition~\ref{def:p-huber-smoothing}). 
This approximation is quadratic near zero, preserves the derivative of the original loss outside this region, and has a uniformly bounded approximation error.

\subsection{Contributions}

Our contributions are as follows.
First, we establish an optimality property of the \(p\)-Huber smoothing construction: for a prescribed smoothness level, it minimizes the uniform approximation error among smooth lower approximations of the scalar loss \(\abs{r}^p/p\). We also give a dual interpretation of this property through its conjugate representation.

Second, we analyze accelerated proximal gradient (APG) with \(p\)-Huber smoothing for proximable convex regularizers and \(1\leq p<2\), including fixed, variable, and backtracking variants. Fixed smoothing attains an \(\varepsilon\)-optimal solution in \(\mathcal O(\varepsilon^{-1/p})\) matrix--vector products. We also establish a lower complexity bound for deterministic methods in our oracle model. This bound matches our upper bound up to constants depending only on \(p\).

Finally, we evaluate the smoothing methods on synthetic regression, UCI regression datasets, sparse text regression, and box-constrained image deblurring. On synthetic data with a small exponent \(p\) and an ill-conditioned matrix, the variable-smoothing method reaches a lower error than the first-order baselines, and it meets a relative accuracy of \(10^{-6}\) on \CtrlVModerateHits{} of the \CtrlTotalInstances{} instances, more than any other method in the comparison. Across the UCI problems, it meets the accuracy targets at least as often as the baselines, and one variant uses fewer matrix--vector products. In the text experiments, our variants produce the smallest final error, with one of them ahead when the penalty is large. The variable-smoothing method also reaches the target accuracy with fewer matrix--vector products in image deblurring tasks.

\subsection{Related work}
\paragraph{First-order methods for H\"older-smooth optimization.}

For \(p\in(1,2)\), the loss \(f\) has a \((p-1)\)-H\"older-continuous gradient. Nesterov's universal accelerated method attains the optimal rate for this generic H\"older class and serves as our canonical black-box baseline \citep{nesterov2015universal}. More recent distance-adaptive methods achieve optimal anytime H\"older rates without prior knowledge of the smoothness parameters or target accuracy; we compare with the deterministic unconstrained AGDA method of \citet{ren2026accelerated}. The closest structure-aware first-order baseline is Generalized-AGD\( + \) of \citet{diakonikolas2024complementary}, whose general-norm framework admits a direct specialization to powered \(\ell_p\) regression. Its guarantee is expressed in the native \(\ell_p\) geometry, whereas ours uses the Euclidean operator norm and initial distance, so neither guarantee uniformly dominates the other.
\paragraph{Smoothing and alternative geometries.} 

Our analysis follows the classical smoothing-and-acceleration paradigm of \citet{nesterov2005smooth} and \citet{beck2012smoothing}. Homotopy methods also use quadratically smoothed \(\ell_p\) objectives \citep{bubeck2018homotopy}. Iterative refinement reduces \(\ell_p\) regression to a short sequence of residual subproblems \citep{adil2019iterative}, and algorithms for \(p\)-norm flows use smoothed \(q\)-norm problems as subroutines \citep{adil2020faster}. Huber-like and piecewise-power losses have also been studied in the context of sparsification \citep{jambulapati2024sparsifying}. A different approach to non-Lipschitz gradients is to change the optimization geometry through relative smoothness and Bregman methods \citep{bauschke2017descent,lu2018relatively}; in particular, \citet{maddison2021dual} develop nonlinear dual-space preconditioning for a complementary \(p\geq2\) regression regime. 

\paragraph{Other high-accuracy methods for standard \(\ell_p\) regression.}\label{sec:related-work-irls}

A substantial literature develops IRLS, homotopy, and iterative-refinement methods with strong complexity guarantees for standard \(\ell_p\) regression \citep{li1993globally,bubeck2018homotopy,adil2019iterative,adil2019fast,adil2024fast,jambulapati2022improved}. Related advances use Lewis weights and sparsification \citep{jambulapati2022improved,jambulapati2024sparsifying}; \citet{ene2026improved} develop a primal--dual IRLS method and bound the number of linear-system solves it requires. They also extend the method to \(1<p<2\) through duality. The quadratically smoothed objectives of \citet{bubeck2018homotopy} and mixed quadratic--\(\ell_p\) residual models in iterative refinement are also closely related to the \(p\)-Huber smoothing used here. These guarantees are typically formulated for unregularized problems and rely on repeated weighted least-squares or linear-system solves. This computational regime is distinct from the composite matrix--vector oracle model considered here.

\section{Problem formulation and \texorpdfstring{\(p\)}{p}-Huber smoothing}\label{sec:problem-setting}

We consider the composite \(\ell_p\) regression problem in~\eqref{eq:lp-regression}, with \(1\leq p<2\), where \(\omega\) is a convex regularizer. We assume that \(\omega:\R^n\to(-\infty,\infty]\) is proper and lower semicontinuous, with an exactly computable proximal operator \(\prox_{\omega/\gamma}\) for every \(\gamma>0\). This includes \(\ell_1\) penalties and indicator functions of simple convex sets. Throughout, \(\norm{\cdot}\) denotes the Euclidean norm on vectors and the spectral norm on matrices.
We assume that the solution set \(\mathcal{U}^\star:=\argmin_{\vu\in\R^n}F(\vu)\) is nonempty and write \(F^\star:=\min_{\vu\in\R^n}F(\vu)\). We also assume that \(\norm{\mA}>0\). Given an initial point \(\vu_0\in\dom(\omega)\), the algorithm is supplied with a bound \(D>0\) satisfying \(D\geq\dist(\vu_0,\mathcal{U}^\star)\).

The computational model restricts access to \(\mA\) to matrix--vector products \(\mA\vu\) and \(\mA^\top\vv\), and provides access to \(\omega\) through its proximal operator. Matrix factorizations and linear-system solves are not treated as unit-cost primitives. Individual matrix--vector products and proximal evaluations are counted separately. The full deterministic oracle model is specified in Subsection~\ref{subsec:lower-bound}.

\subsection{Optimal scalar smoothing}\label{sec:optimal-scalar-smoothing}
The goal is to construct a smoothing surrogate for the composite objective \(F(\vu)=f(\vu)+\omega(\vu)\) in~\eqref{eq:lp-regression}.
We smooth only \(f\) and leave \(\omega\) unchanged. Since \(f\) is separable in the residual coordinates \(r_i=(\mA\vu-\vb)_i\), it is reasonable to first study the scalar loss \(\psi(r):=\abs{r}^p/p\).

For \(L>0\), define the class of \(L\)-smooth functions lying below \(\psi\) by
\[
    \Phi_L := \{\phi\in \mathcal{C}^1(\R): \Lip(\nabla\phi)\leq L, \phi\leq \psi\},
\]
where \(\mathcal{C}^1(\R)\) denotes the class of continuously differentiable functions on \(\R\), and \(\Lip(\nabla\phi)\) denotes the Lipschitz constant of \(\nabla\phi\).
For \(1<p<2\), the curvature \(\nabla^2\psi(r)=(p-1)\abs{r}^{p-2}\) is unbounded near zero but decreases as \(\abs{r}\) increases. This suggests replacing the loss near zero by a quadratic function while preserving its derivative outside this region. Following the smoothing constructions of \citet{adil2019iterative, adil2024fast, bubeck2018homotopy}, we consider the following \(p\)-Huber smoothing of \(\psi\).

\begin{definition}[\(p\)-Huber smoothing]\label{def:p-huber-smoothing}
For \(1\leq p<2\) and a smoothing radius \(\rho>0\), the \(p\)-Huber smoothing surrogate \(\phi_\rho\) of \(\psi(r) = \abs{r}^p/p\) is defined as
\begin{equation}\label{eq:p-huber-smoothing-surrogate}
    \phi_\rho(r) :=
    \begin{cases}
        \frac{1}{2}\rho^{p-2}r^2 & \abs{r} \leq \rho, \\[0.5em]
        \frac{1}{p}\abs{r}^p - c_p\rho^p & \abs{r} > \rho,
    \end{cases} \qquad c_p:=(2-p)/(2p).
\end{equation}
\end{definition} 
For \(\rho>0\), \(\phi_\rho\) is convex and continuously differentiable. Its derivative \(\nabla\phi_\rho(r)=(\max\{\rho,\abs r\})^{p-2}r\) has Lipschitz constant \(\rho^{p-2}\), and the approximation error satisfies \(0\leq\psi-\phi_\rho\leq c_p\rho^p\).
By contrast, for \(1<p<2\) the quadratically penalized Moreau envelope does not keep a uniformly bounded gap.
\begin{proposition}\label{prop:moreau-gap}
For \(L>0\), define the quadratically penalized Moreau envelope of \(\psi\) by
\[
    \tilde\psi_L(r)
    :=
    \inf_{z\in\R}
    \left\{
        \psi(z)+\frac{L}{2}(r-z)^2
    \right\}.
\]
Then \(\tilde\psi_L\leq\psi\). If \(p=1\), then \(\sup_r\{\psi(r)-\tilde\psi_L(r)\}=1/(2L)\). If \(1<p<2\), then \(\sup_r\{\psi(r)-\tilde\psi_L(r)\}=\infty\).
\end{proposition}
For \(1<p<2\), let \(q:=p/(p-1)>2\). The convex conjugate of
\(\psi\) is \(\psi^*(v)=\abs v^q/q\). At \(p=1\),
\(\psi^*=\iota_{[-1,1]}\), where \(\iota_C\) is zero on \(C\)
and \(+\infty\) outside \(C\). The next proposition gives the
conjugate \(\phi_\rho^*\) and its strong convexity and
approximation properties.

\begin{proposition}\label{prop:p-huber-conjugate}
For \(1\leq p<2\) and every \(\rho>0\), the conjugate of
\(\phi_\rho\) has the following form. If \(1<p<2\), then
\begin{equation}\label{eq:p-huber-conjugate}
    \phi_\rho^*(v)=
    \begin{cases}
        \frac{1}{2}\rho^{2-p}v^2
        & \abs v\leq\rho^{p-1}, \\[0.5em]
        \frac{1}{q}\abs v^q+c_p\rho^p
        & \abs v>\rho^{p-1}.
    \end{cases}
\end{equation}
At \(p=1\),
\(\phi_\rho^*(v)=\frac{\rho}{2}v^2+\iota_{[-1,1]}(v)\).
In both cases, \(\phi_\rho^*\) is
\(\rho^{2-p}\)-strongly convex and satisfies
\begin{equation}\label{eq:bounded-dual-correction}
    \psi^*(v)\leq\phi_\rho^*(v)
    \leq\psi^*(v)+c_p\rho^p
    \qquad\text{for all }v\in\R.
\end{equation}
Moreover,
\begin{equation}\label{eq:scalar-dual-p-huber}
    \phi_\rho(r)
    =\max_{v\in\R}\{rv-\phi_\rho^*(v)\}.
\end{equation}
\end{proposition}

By conjugacy, the strong convexity of \(\phi_\rho^*\) gives
the \(\rho^{p-2}\)-smoothness of \(\phi_\rho\), while
\eqref{eq:bounded-dual-correction} gives
\(\psi-c_p\rho^p\leq\phi_\rho\leq\psi\).

For \(L>0\), define the optimal uniform approximation error over
\(\Phi_L\) by
\begin{equation*}
    \delta^\star(L) := \inf_{\phi\in\Phi_L} \sup_{r\in\R} \left\{\psi(r) - \phi(r)\right\}.
\end{equation*}
Following the analytical framework of \citet{samakhoana2026optimal}, the next theorem shows that \(p\)-Huber smoothing attains the minimum uniform approximation error over \(\Phi_L\).

\begin{theorem}\label{thm:optimal-smoothing-surrogate}
    For \(1\leq p<2\) and any \(L>0\), the \(p\)-Huber smoothing surrogate \(\phi_\rho\) defined in~\eqref{eq:p-huber-smoothing-surrogate} with \(\rho = L^{-\frac{1}{2-p}}\)
    attains the optimal uniform error \(\delta^\star(L) = c_p L^{ - \frac{p}{2-p}}\) among all differentiable surrogates in \(\Phi_L\).
\end{theorem}

\subsection{Composite surrogate}

Applying the scalar surrogate to each residual coordinate gives
the smoothed loss
\begin{equation*}
    f_\rho(\vu)
    :=\sum_{i=1}^m\phi_\rho\bigl((\mA\vu-\vb)_i\bigr),
    \qquad \rho>0.
\end{equation*}
The scalar conjugate representation
\eqref{eq:scalar-dual-p-huber} gives
\begin{equation*}
    f_\rho(\vu)
    =\max_{\vv\in\R^m}
    \left\{
        \langle\mA\vu-\vb,\vv\rangle
        -\sum_{i=1}^m\phi_\rho^*(v_i)
    \right\}.
\end{equation*}
The maximizer \(\vv_\rho(\vu)\) is unique and satisfies
\[
    \bigl(\vv_\rho(\vu)\bigr)_i
    =\nabla\phi_\rho\bigl((\mA\vu-\vb)_i\bigr),
    \qquad i=1,\ldots,m.
\]
Hence,
\(\nabla f_\rho(\vu)=\mA^\top\vv_\rho(\vu)\).

\begin{proposition}\label{prop:p-huber-properties}
    The function \(f_\rho\) is convex and continuously differentiable, with \(\Lip(\nabla f_\rho) \leq \norm{\mA}^2 \rho^{p-2}\).
    The approximation error satisfies \(0 \leq f(\vu) - f_\rho(\vu) \leq c_p m \rho^p\) for all \(\vu\in\R^n\).
    Moreover, the surrogate \(f_\rho\) is non-increasing in \(\rho\), i.e., for any \(0 < \rho_2 \leq \rho_1\) and any \(\vu \in \R^n\),
    \begin{equation}\label{eq:p-huber-radius-drift}
        0 \leq f_{\rho_2}(\vu) - f_{\rho_1}(\vu) \leq c_p m (\rho_1^p - \rho_2^p).
    \end{equation}
\end{proposition}

Adding \(\omega\) gives the composite surrogate
\begin{equation*}
    F_\rho(\vu) := f_\rho(\vu) + \omega(\vu).
\end{equation*}

\begin{proposition}\label{prop:composite-surrogate-properties}
For every \(\rho>0\), the function \(F_\rho\) is proper, lower
semicontinuous, and convex.  For every \(\vu\in\dom(\omega)\),
\(0\leq F(\vu)-F_\rho(\vu)\leq c_p m\rho^p\).
For any \(0<\rho_2\leq\rho_1\) and
\(\vu\in\dom(\omega)\),
\(0\leq F_{\rho_2}(\vu)-F_{\rho_1}(\vu)
\leq c_p m\bigl(\rho_1^p-\rho_2^p\bigr)\).
Moreover, for every \(\vu\in\dom(\omega)\) and
\(\vu^\star\in\mathcal{U}^\star\),
\(F(\vu)-F^\star\leq c_p m\rho^p
+F_\rho(\vu)-F_\rho(\vu^\star)\).
\end{proposition}

\section{Optimal deterministic matrix--vector complexity}\label{sec:optimal-matrix-vector-complexity}

We apply accelerated proximal gradient to the \(p\)-Huber surrogate and analyze its convergence under fixed and variable smoothing schedules.
Subsection~\ref{subsec:hs-apg} specifies the algorithm, and Subsection~\ref{subsec:upper-bound} analyzes the fixed-smoothing and variable-smoothing variants. Subsection~\ref{subsec:lower-bound} establishes a lower complexity bound. Together with the upper bound, it determines the optimal deterministic matrix--vector complexity of \(\ell_p\) regression for \(1\leq p<2\).

\subsection{Accelerated proximal gradient with \texorpdfstring{\(p\)}{p}-Huber smoothing}
\label{subsec:hs-apg}

Algorithm~\ref{alg:hs-apg} applies accelerated proximal gradient to the \(p\)-Huber surrogate using Nesterov's three-sequence scheme \citep{nesterov2005smooth} and a non-increasing sequence of positive smoothing radii \(\{\rho_t\}_{t=0}^T\). Write \(f_t:=f_{\rho_t}\), \(F_t:=F_{\rho_t}\), and \(L_t:=\norm{\mA}^2\rho_t^{p-2}\). For the upper-bound algorithms and proofs below, a known upper bound may be used in place of \(\norm{\mA}\), with the corresponding guarantees expressed in terms of that bound. For \(t = 1,\ldots,T\), define \(\Gamma_t := t(t+1)/2\) and \(\beta_t := \alpha_t\gamma_t\Gamma_t\), where \(\alpha_t\) and \(\gamma_t\) are the coefficients in Algorithm~\ref{alg:hs-apg}, and set \(\Gamma_0=0\). The following lemma gives a one-step bound used in the convergence analysis of both fixed and variable smoothing.

\begin{algorithm}[t]
\caption{Accelerated proximal gradient with \(p\)-Huber smoothing}
\label{alg:hs-apg}
\KwIn{forward/adjoint access to \(\mA\), \(\vb\), \(p\in[1,2)\), initial point \(\vu_0\), positive smoothing radii \(\{\rho_t\}_{t=0}^T\), and positive parameters \(\{\gamma_t\}_{t=1}^T\)}
\KwInit{\(\vy_0 \leftarrow \vu_0\), \(\vz_0 \leftarrow \vu_0\)}
\For{\(t=1,\ldots,T\)}{
    \(\alpha_t \leftarrow 2/(t+1)\)\;
    \(\vx_t \leftarrow (1-\alpha_t)\vy_{t-1}+\alpha_t\vz_{t-1}\)\;
    \(\vg_t \leftarrow \nabla f_{\rho_t}(\vx_t) = \mA^\top\nabla\phi_{\rho_t}(\mA\vx_t-\vb)\)\;
    \(\vz_t \leftarrow \prox_{\omega/\gamma_t}(\vz_{t-1}-\gamma_t^{-1}\vg_t)\)\;
    \(\vy_t \leftarrow (1-\alpha_t)\vy_{t-1}+\alpha_t\vz_t\)\;
}
\KwOut{\(\vy_T\)}
\end{algorithm}

\begin{lemma}\label{lem:accelerated-drift}
Suppose that \(\gamma_t\geq\alpha_t L_t\) and \(\rho_t\leq\rho_{t-1}\) for all \(t=1,\ldots,T\). Then, for every \(\vu\in\dom(\omega)\) and \(t=1,\ldots,T\), Algorithm~\ref{alg:hs-apg} satisfies
\begin{align}
\label{eq:accelerated-drift}
&\Gamma_t\bigl(F_t(\vy_t)-F_t(\vu)\bigr)
~\leq~
\Gamma_{t-1}\bigl(F_{t-1}(\vy_{t-1})-
                         F_{t-1}(\vu)\bigr)
\\ & \qquad +\frac{\beta_t}{2}
 \left(\norm{\vu-\vz_{t-1}}^2-\norm{\vu-\vz_t}^2\right)
+ c_p m\Gamma_{t-1}\bigl(\rho_{t-1}^p-\rho_t^p\bigr).
\nonumber
\end{align}
\end{lemma}

The proof combines the smooth upper model for \(f_t\), convexity of \(f_t\) and \(\omega\), the proximal three-point inequality, and the bound on \(F_t-F_{t-1}\) in Proposition~\ref{prop:composite-surrogate-properties}; see Appendix~\ref{app:proof-accelerated-drift}.

\subsection{Upper bound on matrix--vector complexity}\label{subsec:upper-bound}

The upper-bound analysis considers fixed and variable \(p\)-Huber smoothing. Fixed smoothing chooses its radius for a prescribed terminal horizon \(T\). Variable smoothing uses a decreasing radius sequence independent of the terminal horizon, with a convergence guarantee under a bounded-iterate assumption.

\paragraph{Fixed smoothing.}
The fixed-smoothing variant optimizes a single surrogate.  The step-size choice below makes
\(\beta_t\) constant, so the distance terms in Lemma~\ref{lem:accelerated-drift} telescope exactly. 

\begin{theorem}
\label{thm:upper-bound}
Let \(1\leq p<2\) and \(T\geq1\), suppose \(\norm{\mA}>0\), and let \(D>0\) satisfy \(D\geq\dist(\vu_0,\mathcal{U}^\star)\).  Set
\begin{equation}
\label{eq:fixed-radius}
    \rho_T
    :=
    \frac{2\norm{\mA}D}{\sqrt{mT(T+1)}},
    \qquad
    L_T:=\norm{\mA}^2\rho_T^{p-2},
\end{equation}
and run Algorithm~\ref{alg:hs-apg} with
\(\rho_t\equiv\rho_T\) and \(\gamma_t:=\frac{2L_T}{t}\), for \(t = 1,\ldots,T\).
Then the output \(\vy_T\) satisfies
\begin{equation}
\label{eq:fixed-rate}
    F(\vy_T)-F^\star
    \leq
    \frac{2^p}{p}\,
    m^{\frac{2-p}{2}}
    \left(
        \frac{\norm{\mA}D}{\sqrt{T(T+1)}}
    \right)^p
    \leq
    \frac{2^p}{p}\,
    m^{\frac{2-p}{2}}
    \left(\frac{\norm{\mA}D}{T}\right)^p.
\end{equation}
\end{theorem}

The rate reflects the optimal scalar error--smoothness trade-off: the optimization and smoothing error bounds are of orders \(\norm{\mA}^2D^2\rho^{p-2}/T^2\) and \(m\rho^p\), respectively. Balancing these terms yields the rate \(\mathcal{O}(T^{-p})\). Prior knowledge of the terminal horizon can be removed by a restart-doubling scheme without changing the order of the rate; details are given in Appendix~\ref{app:variants}.

\paragraph{Variable smoothing.}
The variable-smoothing variant uses the radius schedule \(\rho_t:=\bar\rho/(t+1)\), which does not require prior knowledge of the total iteration count \(T\). Since \(L_t=\norm{\mA}^2\bar\rho^{p-2}(t+1)^{2-p}\), the upper bound on the gradient Lipschitz constant grows as \((t+1)^{2-p}\), while the uniform approximation error bound decreases as \((t+1)^{-p}\).

\begin{theorem}\label{thm:variable-radius}
Let \(1\leq p<2\) and \(T \geq 1\), and suppose \(\norm{\mA}>0\). Suppose further that for some minimizer \(\vu^\star\in\mathcal{U}^\star\) and some \(R>0\), the accelerated sequence satisfies \(\norm{\vz_t-\vu^\star}\leq R\) for all \(t=0,\ldots,T\).
Set
\begin{equation}
\label{eq:variable-rbar}
    \bar\rho
    :=
    2\sqrt{\frac{2-p}{(p+2)m}}\,\norm{\mA}R,
    \qquad
    \rho_t:=\frac{\bar\rho}{t+1},
    \qquad
    L_t:=\norm{\mA}^2\rho_t^{p-2},
\end{equation}
and run Algorithm~\ref{alg:hs-apg} with \(\gamma_t:=\alpha_t L_t\).
Then, after \(T\) iterations,
\begin{equation}
\label{eq:variable-rate}
    F(\vy_T)-F^\star
    \leq
    C_p\,
    m^{\frac{2-p}{2}}
    \left(\frac{\norm{\mA}R}{T+1}\right)^p,
    \qquad
    C_p
    :=
    \frac{4}{p}
    \left(\frac{p+2}{4(2-p)}\right)^{\frac{2-p}{2}}.
\end{equation}
\end{theorem}

For either variant, \(T\) iterations of Algorithm~\ref{alg:hs-apg} use \(2T\) individual matrix--vector products and \(T\) proximal evaluations.

\paragraph{Practical variants.}
In practice, there are several heuristic ways to further improve the performance.
First, the backtracking variant replaces the global bound \(L_t\) by a local estimate \(\widehat L_t\), accepting a candidate whenever
\begin{equation}
\label{eq:p-huber-line-search-test}
f_t(\vy_t)
\leq
f_t(\vx_t)
+\ip{\nabla f_t(\vx_t)}{\vy_t-\vx_t}
+\frac{\widehat L_t}{2}\norm{\vy_t-\vx_t}^2.
\end{equation}
Appendix~\ref{app:backtracking-line-search} gives additional conditions under which this variant retains the variable-smoothing convergence bound.

Second, for problem~\eqref{eq:lp-regression} with \(\omega\equiv0\), we also consider initialization at a least-squares solution \(\vu_{\mathrm{OLS}}\in\arg\min_{\vu}\frac12\|\mA\vu-\vb\|_2^2\). This allows the method to initialize without prior knowledge of the distance \(D\). If \(\mA\vu_{\mathrm{OLS}}=\vb\), the problem is solved. Otherwise, the variable-smoothing variant starts at \(\vu_{\mathrm{OLS}}\) and sets \(\bar\rho=2\|\mA\vu_{\mathrm{OLS}}-\vb\|_\infty\). It uses \(\rho_t=\bar\rho/(t+1)\) and the global bound \(L_t\), without backtracking or reference-distance information.  The full specification and its connection to the variable-smoothing bound are given in Appendix~\ref{app:ols-variant}.

\subsection{Lower bound on matrix--vector complexity}\label{subsec:lower-bound}

The lower-bound analysis follows the framework of \citet{ouyang2021lower} and applies to all adaptive deterministic methods in the following oracle model. 

The algorithm accesses \(\mA\) through the paired matrix--vector oracle
\[
\oracle_{\mA}(\vu,\vv) := (\mA\vu,\mA^\top\vv), \qquad \vu\in\R^n,\quad \vv\in\R^m.
\]
For a query budget \(Q\in\mathbb{N}\), the algorithm is given \((p,m,n,\vb,\vu_0,D,\norm{\mA},Q)\) before its first query. Each query and the final output may depend on these data, all preceding oracle responses, and all proximal evaluations made so far.
Between queries, the algorithm may perform arbitrary deterministic computations and any number of exact evaluations of \(\prox_{\omega/\gamma}\), for any \(\gamma>0\). Apart from the supplied norm \(\norm{\mA}\), information about \(\mA\) is available only through the oracle. Matrix inversion, factorization, and linear-system solution are not available as unit-cost primitives.

To establish a worst-case lower bound for the composite problem class, it suffices to consider instances with \(\omega\equiv0\) in Theorem~\ref{thm:matvec-lower-bound}; the proof is given in Appendix~\ref{app:proof-matvec-lower-bound}.

\begin{theorem}\label{thm:matvec-lower-bound}
Fix \(p\in[1,2)\), \(\Lambda,D>0\), and an integer \(Q\geq1\). Suppose \(m\geq8Q+18\) and \(n\geq4Q+10\). Consider any adaptive deterministic algorithm initialized at \(\vu_0\in\R^n\) and making at most \(Q\) calls to the paired oracle, with any number of exact proximal evaluations. There exist \(\mA\in\R^{m\times n}\) and \(\vb\in\R^m\) for problem~\eqref{eq:lp-regression} with \(\omega\equiv0\) such that \(\norm{\mA}=\Lambda\), \(\dist(\vu_0,\mathcal{U}^\star)=D\), and \(F^\star=0\). On this instance, the algorithm's output \(\vu_Q\) satisfies
\begin{equation}
\label{eq:matvec-lower-bound}
    F(\vu_Q)-F^\star
    \geq
    \frac{1}{p}
    \left(
        \frac{\sqrt{3}\,c_0}{4\sqrt{2}}
    \right)^{\!p}
    m^{\frac{2-p}{2}}
    \left(
        \frac{\Lambda D}{2Q+5}
    \right)^{\!p},
\end{equation}
where \(c_0\in(0,1)\) is a universal constant.
\end{theorem}

For an integer \(T\geq1\), under the theorem's assumptions with \(Q=2T\), this gives a worst-case objective gap of \(\Omega(T^{-p})\) for any method using at most \(2T\) individual products. Since Algorithm~\ref{alg:hs-apg} with fixed smoothing attains the upper bound in Theorem~\ref{thm:upper-bound} with \(2T\) products, its deterministic matrix--vector complexity is optimal up to constants depending only on \(p\) in the stated oracle model.

\section{Numerical experiments}\label{sec:experiments}

We compare methods based on \(p\)-Huber smoothing with first-order baselines using objective accuracy and matrix--vector product counts. Subsection~\ref{subsec:unregularized-regression} evaluates regression without regularization (i.e., $\omega \equiv 0$) on synthetic and real-world data. Subsection~\ref{subsec:composite-regularization} evaluates regression with \(\ell_1\) regularization.

Accuracy is measured by the relative objective gap
\( \operatorname{RelGap}(\vu) = (F(\vu)-F_{\mathrm{ref}})/(F(\vu_0)-F_{\mathrm{ref}}) \),
where \(F_{\mathrm{ref}}\) is a reference objective value and \(\vu_0\) is a common normalization point for all methods on each problem. Computational cost is measured by the number of matrix--vector products, including those used for initialization and line search. Reference construction and implementation details are provided in Appendix~\ref{app:experiment-details}.

\subsection{Regression without regularization}\label{subsec:unregularized-regression}

\subsubsection{Synthetic data}\label{subsec:controlled-studies}

We study how \(p\) and the conditioning of \(\mA\) affect the performance of fixed smoothing, variable smoothing, and variable smoothing with backtracking. We compare these methods with Nesterov's Universal Fast Gradient Method (UFGM) \citep{nesterov2015universal}, the Accelerated Distance-Adaptive method (AGDA) \citep{ren2026accelerated}, and Generalized-AGD+ (GAGD+) \citep{diakonikolas2024complementary}.

The synthetic instances cover all combinations of \(p\in\{1.1,1.25,1.5,\)\allowbreak\(1.75,1.9\}\) and \(\kappa(\mA)\in\{1,10,\)\allowbreak\(10^2,10^3,10^4\}\). We generate ten instances for each pair, giving 250 instances in total, with \(m=200\), \(n=50\), and \(\norm{\mA}=1\). All methods start from \(\vu_0=\mathbf 0\), with initial distance \(\norm{\vu_0-\vu^\star}=D=10\). The construction provides the exact reference value \(F_{\mathrm{ref}}=F^\star\) for each instance; details are given in Appendix~\ref{app:controlled-instance-details}.
Each method receives a maximum budget of \(60{,}000\) matrix--vector products. Fixed smoothing uses \(T=30{,}000\) to select its smoothing radius, while UFGM and GAGD+ are configured for a relative-gap target of \(10^{-9}\). For each method and instance, we evaluate attainment of \(10^{-6}\) and \(10^{-9}\) along a single run. Further algorithm settings are given in Appendix~\ref{app:method-details}. In the figures and tables, HS-APG-F, HS-APG-V, and HS-APG-V-BT denote fixed smoothing, variable smoothing, and variable smoothing with backtracking. HS-APG-OLS denotes variable smoothing with least-squares initialization and is used in the real-data experiments.

\begin{figure*}[t]
    \centering
    \includegraphics[width=\textwidth]{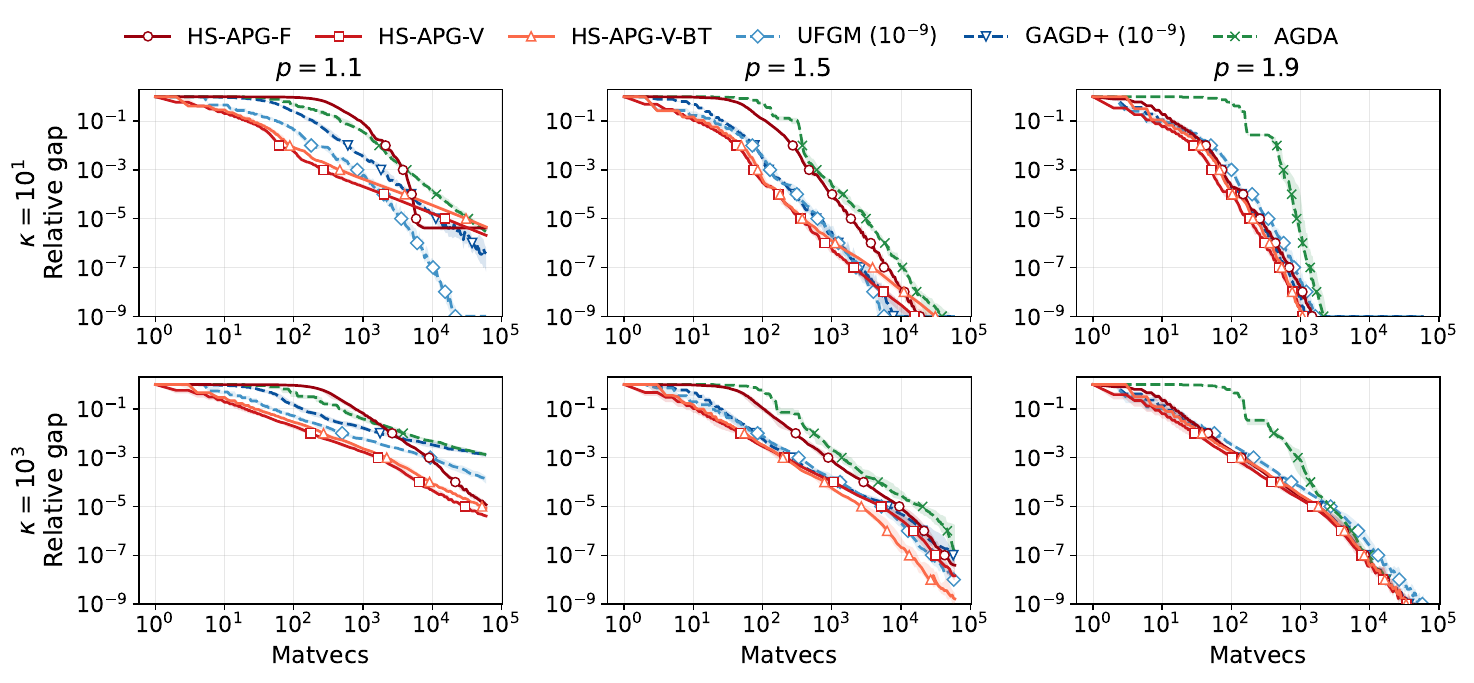}
    \caption{Convergence on synthetic instances with
    \(p\in\{1.1,1.5,1.9\}\) and
    \(\kappa(\mA)\in\{10,10^3\}\). Curves show median relative gaps over ten
    paired seeds against the number of matrix--vector products. Shaded
    bands span the 25th to 75th percentiles. UFGM and GAGD+ use target
    \(10^{-9}\).
    The complete \(5\times5\)
    grid is shown in Figure~\ref{fig:controlled-full-grid} of the appendix.}
    \label{fig:controlled-convergence}
\end{figure*}

Figure~\ref{fig:controlled-convergence} compares convergence for selected exponents and condition numbers; full-grid results are provided in Appendix~\ref{app:results-details}. Variable smoothing achieves consistently better terminal accuracy than the baselines in the small-\(p\), poorly conditioned regime, yielding lower terminal gaps on all \CtrlLowPHighKCount{} instances with \(p\leq1.25\) and \(\kappa(\mA)\geq10^2\). Across the full grid, it also reaches \(10^{-6}\) on \CtrlVModerateHits{} of \CtrlTotalInstances{} instances, the highest attainment count among the six methods. UFGM achieves lower median terminal gaps at \(p=1.1\) and \(\kappa(\mA)\in\{1,10\}\), and leads at the stricter target \(10^{-9}\), reaching \CtrlUFGMStrictHits{} of \CtrlTotalInstances{} instances. At \(\kappa(\mA)=10^4\), high accuracy remains difficult for all methods: each reaches \(10^{-9}\) on at most two of the 50 instances (Table~\ref{tab:controlled-attainment}).

The variable-smoothing methods generally achieve lower terminal gaps than fixed smoothing, but neither variant uniformly dominates the other. Without backtracking, variable smoothing has lower terminal gaps on all \(p=1.1\) instances, whereas backtracking achieves a lower median terminal gap at \(p=1.5\) and \(\kappa(\mA)=10^3\). Backtracking can permit larger step sizes, but its line search requires additional matrix--vector products within the same budget.

\subsubsection{Real-world data}
\label{subsec:real-data}

We evaluate the methods on six UCI regression datasets (Power, Year,
Protein, KEGG, CT-Slices, and Buzz) with \(p\in\{1.1,1.5,1.9\}\), yielding
18 optimization problems. Since \(D\) may be unknown in practice, we also
test a least-squares-initialized variant of variable smoothing. It sets
the smoothing scale from the least-squares residual and requires no
reference-distance information.
The least-squares solver starts at zero, as do all other methods. The other smoothing variants use
\(D_{\mathrm{ref}}=\|\vu_{\mathrm{ref}}-\vu_0\|_2\), the distance from the
initial point to a numerical reference candidate, to set their radii.

Each method receives a budget of \(60{,}000\) matrix--vector products,
including initialization and line search. Fixed smoothing uses
\(T=30{,}000\) to set its radius. We compare attainment at relative-gap
targets \(10^{-4}\) and \(10^{-6}\). UFGM and GAGD+ are configured for
the accuracy target used in each comparison. Dataset sources,
preprocessing, reference certification, and implementation details are
given in Appendix~\ref{app:realdata-experiments}.

\begin{table}[t]
    \caption{Target attainment on real data within a budget of 60,000 algorithmic matrix--vector products. Hit counts cover all 18 problems. Median first-hit costs use the same subset for every method: the 15 problems at $10^{-4}$ and 13 at $10^{-6}$ where all seven methods reach the target. UFGM and GAGD+ use target-matched configurations. Costs include initialization and line search; products used only for objective reporting are excluded for every method.}
    \label{tab:realdata-summary}
    \centering
    \setlength{\tabcolsep}{9pt}
    \begin{tabular}{lrrrr}
        \toprule
        & \multicolumn{2}{c}{Target $10^{-4}$} & \multicolumn{2}{c}{Target $10^{-6}$}\\
        \cmidrule(lr){2-3}\cmidrule(lr){4-5}
        Method & Hits / 18 & Median matvecs & Hits / 18 & Median matvecs\\
        \midrule
        HS-APG-F & 17 & 2,498 & 15 & 5,900\\
        HS-APG-V & 17 & 576 & 15 & 1,556\\
        HS-APG-V-BT & \textbf{18} & 675 & 15 & 875\\
        HS-APG-OLS & 17 & \textbf{241} & \textbf{16} & \textbf{673}\\
        \midrule
        UFGM & 17 & 978 & 14 & 1,167\\
        GAGD+ & 16 & 738 & 13 & 1,734\\
        AGDA & 16 & 1,564 & 13 & 2,241\\
        \bottomrule
    \end{tabular}
\end{table}

Table~\ref{tab:realdata-summary} shows that the smoothing methods reach each target on at least as many problems as the baselines, while variable smoothing has lower median costs than fixed smoothing. Backtracking reaches \(10^{-4}\) on all \RDProblemCount{} problems. At \(10^{-6}\), each smoothing variant reaches the target on at least \RDBTStrictHits{} problems, compared with \RDGAGDStrictHits{} or \RDUFGMStrictHits{} for the baselines. Both variable-smoothing variants have lower median costs than all three baselines at \(10^{-4}\); for the backtracking variant, this advantage also holds at \(10^{-6}\). These cost comparisons use the common subsets where all seven methods reach the corresponding target. The least-squares-initialized variant reaches the two targets on \RDOLSLooseHits{} and \RDOLSStrictHits{} problems, respectively, with low median costs including initialization. It provides a practical option when \(D\) is unknown. Tables~\ref{tab:realdata-refdist-attainment-1e-4} and~\ref{tab:realdata-refdist-attainment-1e-6} provide problem-level results, while Appendix~\ref{app:realdata-trajectories} presents convergence curves and a separate gradient-restart ablation.

\subsection{Regression with \texorpdfstring{\(\ell_1\)}{l1} regularization}\label{subsec:composite-regularization}

\newcommand{\EThreeCompleted}{48}
\newcommand{\EThreeHitsFour}{0}
\newcommand{\EThreeHitsSix}{0}

We evaluate the composite formulation on \(\ell_1\)-regularized regression using E2006 financial-report text \citep{e2006textregression}. The smoothing methods apply the \(p\)-Huber surrogate to the loss and handle the penalty \(\omega(\vu)=\lambda\norm{\vw}_1\) on the feature coefficients \(\vw\) by soft thresholding. We use the TFIDF and log1p representations with \(p=1.5\), an unpenalized intercept, and \(\lambda/\lambda_{\max}\in\{0.01,0.03,0.1,0.3\}\), giving eight optimization problems. Here, \(\lambda_{\max}\) is the smallest penalty for which setting all feature coefficients to zero is optimal \citep{o2016statistical}; its precise definition is given in~\eqref{eq:e2006-lambda-max} in Appendix~\ref{app:regularized_setup}. We omit HS-APG-F in this experiment because the preceding studies already evaluate the fixed-smoothing method extensively and show that the variable-smoothing variants are generally more competitive in first-hit cost. The E2006 study therefore focuses on the practical variants HS-APG-V and HS-APG-V-BT, together with composite UFGM and GAGD+. The two baselines are each run separately for relative-gap targets \(10^{-4}\) and \(10^{-6}\).

All methods start with zero feature coefficients and an intercept chosen to minimize the loss. Each run has a budget of \(30{,}000\) matrix--vector products for TFIDF or \(12{,}000\) for log1p, including line-search trials. The smoothing variants set their smoothing radii using the distance to a reference candidate as a proxy, and relative gaps are computed using independently certified reference lower bounds. Appendix~\ref{app:regularized_setup} gives the preprocessing, reference construction, and distance proxy.

Table~\ref{tab:e2006-grid-terminal} reports the terminal relative gap of every run. Variable smoothing attains the smallest terminal gap on TFIDF at ratios \(0.01\), \(0.03\), and \(0.1\), and on log1p at \(0.01\) and \(0.03\). Variable smoothing with backtracking attains the smallest terminal gap in the remaining cases: \(0.3\) on TFIDF, and \(0.1\) and \(0.3\) on log1p. The benefit of backtracking over variable smoothing is therefore concentrated at the larger penalties in this grid. No run reaches a relative gap of \(10^{-4}\) within its budget. Appendix~\ref{app:regularized_setup} reports thresholded active-feature counts. 

Appendix~\ref{app:constrained_setup} complements the \(\ell_1\)-regularized study with an image deblurring experiment using \(\omega=\iota_{[0,1]^n}\).

\begin{table}[H]
\centering
\setlength{\tabcolsep}{6pt}
\caption{Terminal relative gaps for all 48 runs. Columns correspond to penalty ratios $\lambda/\lambda_{\max}$. Budgets are $30{,}000$ matrix--vector products for TFIDF and $12{,}000$ for log1p. Bold entries indicate the smallest gap among the six configurations in each column within each representation.}
\label{tab:e2006-grid-terminal}
\begin{tabular}{lrrrr}
\toprule
\multicolumn{5}{c}{TFIDF} \\
Method / target & $0.01$ & $0.03$ & $0.1$ & $0.3$ \\
\midrule
HS-APG-V & $\mathbf{3.0992\times10^{-3}}$ & $\mathbf{2.0121\times10^{-3}}$ & $\mathbf{7.8028\times10^{-4}}$ & $1.6917\times10^{-3}$ \\
HS-APG-V-BT & $4.9278\times10^{-3}$ & $2.9283\times10^{-3}$ & $1.7572\times10^{-3}$ & $\mathbf{4.6715\times10^{-4}}$ \\
UFGM ($10^{-4}$) & $1.2673\times10^{-2}$ & $8.7801\times10^{-3}$ & $5.9847\times10^{-3}$ & $2.8517\times10^{-3}$ \\
UFGM ($10^{-6}$) & $1.2660\times10^{-2}$ & $8.7173\times10^{-3}$ & $6.0044\times10^{-3}$ & $2.8439\times10^{-3}$ \\
GAGD+ ($10^{-4}$) & $6.7852\times10^{-3}$ & $3.7290\times10^{-3}$ & $9.4428\times10^{-4}$ & $7.1823\times10^{-4}$ \\
GAGD+ ($10^{-6}$) & $6.7852\times10^{-3}$ & $3.7290\times10^{-3}$ & $9.4428\times10^{-4}$ & $7.1823\times10^{-4}$ \\
\midrule
\multicolumn{5}{c}{log1p} \\
Method / target & $0.01$ & $0.03$ & $0.1$ & $0.3$ \\
\midrule
HS-APG-V & $\mathbf{1.1392\times10^{-2}}$ & $\mathbf{3.2231\times10^{-2}}$ & $3.3366\times10^{-2}$ & $8.0258\times10^{-3}$ \\
HS-APG-V-BT & $1.3119\times10^{-2}$ & $3.5383\times10^{-2}$ & $\mathbf{1.2197\times10^{-2}}$ & $\mathbf{1.3372\times10^{-3}}$ \\
UFGM ($10^{-4}$) & $2.1655\times10^{-2}$ & $5.9475\times10^{-2}$ & $8.4516\times10^{-2}$ & $5.6877\times10^{-3}$ \\
UFGM ($10^{-6}$) & $2.1670\times10^{-2}$ & $5.9423\times10^{-2}$ & $8.4810\times10^{-2}$ & $5.6795\times10^{-3}$ \\
GAGD+ ($10^{-4}$) & $1.6539\times10^{-2}$ & $4.4306\times10^{-2}$ & $1.6620\times10^{-2}$ & $1.9741\times10^{-3}$ \\
GAGD+ ($10^{-6}$) & $1.6547\times10^{-2}$ & $4.4306\times10^{-2}$ & $1.6620\times10^{-2}$ & $1.9741\times10^{-3}$ \\
\bottomrule
\end{tabular}
\end{table}

\section{Conclusion}\label{sec:conclusion}

We establish the optimal one-sided approximation--smoothness trade-off for
\(p\)-Huber smoothing and apply the resulting surrogate in APG for composite \(\ell_p\) regression.
The fixed-smoothing variant attains an \(\varepsilon\)-optimal solution for \(1\leq p<2\) in
\(\mathcal O(\varepsilon^{-1/p})\) matrix--vector products, matching the
deterministic lower bound in the stated oracle model.
Experiments on synthetic and real data show
competitive performance and illustrate the method's applicability to
\(\ell_1\) regularization and box constraints. 

\bibliography{reference}
\bibliographystyle{plainnat}

\appendix

\section{Omitted proofs from the main body}\label{app:proofs}

This appendix provides the proofs of the results stated in the main body. 

\subsection{Proofs for Section~\ref{sec:problem-setting}}\label{app:proofs-problem-setting}

\subsubsection{Properties of the smoothed regression loss}
\label{app:proof-coordinatewise-lifting}

The following lemma transfers the smoothness and approximation bounds of
the scalar surrogate to the regression loss. We use it to prove
Proposition~\ref{prop:p-huber-properties}.

\begin{lemma}\label{lem:coordinatewise-lifting}
Let \(L>0\) and let \(\phi\in \mathcal{C}^1(\R)\) be convex with \(\Lip(\nabla\phi)\leq L\). Suppose
\(0\leq\psi(r)-\phi(r)\leq\delta\) for all \(r\in\R\). Then
\[
    f_\phi(\vu):=\sum_{i=1}^m
    \phi\bigl((\mA\vu-\vb)_i\bigr)
\]
is convex and continuously differentiable, with
\(\Lip(\nabla f_\phi)\leq L\norm{\mA}^2\) and
\(0\leq f(\vu)-f_\phi(\vu)\leq m\delta\) for all \(\vu\in\R^n\).
\end{lemma}

\begin{proof}
Convexity follows because each summand is the composition of a convex scalar
function with an affine map.  The chain rule gives
\[
    \nabla f_\phi(\vu)
    =\mA^\top\nabla\phi(\mA\vu-\vb),
\]
where \(\nabla\phi\) is applied coordinatewise.  Thus, for all
\(\vu,\vw\in\R^n\),
\begin{align*}
\norm{\nabla f_\phi(\vu)-\nabla f_\phi(\vw)}
&\leq
\norm{\mA}\,
\norm{\nabla\phi(\mA\vu-\vb)-
      \nabla\phi(\mA\vw-\vb)} \\
&\leq
L\norm{\mA}\,\norm{\mA(\vu-\vw)}
\leq L\norm{\mA}^2\norm{\vu-\vw}.
\end{align*}
Finally, summing
\(0\leq\psi(r_i)-\phi(r_i)\leq\delta\) over the \(m\) residual
coordinates yields the uniform approximation bound.
\end{proof}

\subsubsection{Proof of Proposition~\ref{prop:moreau-gap}}
\label{app:proof-moreau-gap}

\begin{proof}
The envelope \(\tilde\psi_L\) is convex and has an \(L\)-Lipschitz derivative, and \(\tilde\psi_L\leq\psi\) because \(z=r\) is feasible. For \(r>0\) it is enough to minimize over \(z\geq0\); the case \(r<0\) follows by evenness.

If \(p=1\), the minimizer is \(z=0\) when \(0\leq r\leq 1/L\), and \(z=r-1/L\) when \(r>1/L\). In the first range the gap is \(r-(L/2)r^2\leq 1/(2L)\). In the second range it equals
\[
    r-\bigl(r-1/L\bigr)-(L/2)(1/L)^2=1/(2L).
\]
Thus the supremum is \(1/(2L)\).

Now let \(1<p<2\). For each \(z>0\) set \(r=z+z^{p-1}/L\). This \(z\) satisfies the first-order condition for \(\tilde\psi_L(r)\), so
\[
    \psi(r)-\tilde\psi_L(r)
    =
    \frac{r^p-z^p}{p}-\frac{z^{2p-2}}{2L}.
\]
The mean-value theorem gives \(r^p-z^p=p\xi^{p-1}(r-z)\) for some \(\xi\in(z,r)\). Since \(r-z=z^{p-1}/L\) and \(\xi>z\),
\[
    \psi(r)-\tilde\psi_L(r)
    >
    z^{p-1}\cdot\frac{z^{p-1}}{L}-\frac{z^{2p-2}}{2L}
    =
    \frac{z^{2p-2}}{2L}.
\]
The right-hand side tends to infinity as \(z\to\infty\), and the corresponding \(r\) does as well. Therefore the supremum is infinite.
\end{proof}

\subsubsection{Proof of Proposition~\ref{prop:p-huber-conjugate}}
\label{app:proof-p-huber-conjugate}

\begin{proof}
At \(p=1\), the conjugate of \(\psi(r)=\abs r\) is \(\psi^*=\iota_{[-1,1]}\). If \(\abs v>1\), the expression \(vr-\phi_\rho(r)\) tends to \(+\infty\) as \(r\) tends to infinity with the sign of \(v\). If \(\abs v\leq1\), the supremum is attained at \(r=\rho v\) and equals \(\rho v^2/2\). This proves the stated endpoint formula. Subtracting \(\rho v^2/2\) from \(\phi_\rho^*\) leaves the convex indicator \(\iota_{[-1,1]}\), so \(\phi_\rho^*\) is \(\rho\)-strongly convex. On \([-1,1]\), its correction to \(\psi^*\) lies in \([0,\rho/2]\); outside this interval, both functions are \(+\infty\). Thus~\eqref{eq:bounded-dual-correction} holds as an extended-real inequality.

Now let \(1<p<2\). Both \(\phi_\rho\) and its conjugate are even, so it suffices to take
\(v\geq0\). By definition,
\[
    \phi_\rho^*(v)=\sup_{r\geq0}\{vr-\phi_\rho(r)\}.
\]
If the maximizer lies in \(0\leq r\leq\rho\), the first-order condition is
\(v=\rho^{p-2}r\), hence \(r=\rho^{2-p}v\). The value \(r=\rho^{2-p}v\)
lies in \([0,\rho]\) if and only if \(0\leq v\leq\rho^{p-1}\), and substitution gives
\(\phi_\rho^*(v)=\frac12\rho^{2-p}v^2\). If the maximizer lies in \(r>\rho\),
the first-order condition is \(v=r^{p-1}\), hence \(r=v^{q-1}\). The value
\(r=v^{q-1}\) satisfies \(r>\rho\) if and only if \(v>\rho^{p-1}\), and
\[
    \phi_\rho^*(v)
    =v^q-\frac1p v^q+c_p\rho^p
    =\frac1q v^q+c_p\rho^p.
\]
At \(v=\rho^{p-1}\), the two expressions agree because
\(1/q+c_p=1/2\). This proves~\eqref{eq:p-huber-conjugate}.

For \(\abs{v}<\rho^{p-1}\), we have \(\nabla^2 \phi_\rho^*(v)=\rho^{2-p}\).
For \(\abs{v}>\rho^{p-1}\),
\[
    \nabla^2 \phi_\rho^*(v)
    =(q-1)\abs{v}^{q-2}
    \geq(q-1)\rho^{(p-1)(q-2)}
    =(q-1)\rho^{2-p}
    \geq\rho^{2-p}.
\]
Since the first derivatives match at \(v=\pm\rho^{p-1}\), \(\phi_\rho^*\) is
\(\rho^{2-p}\)-strongly convex on \(\R\).

It remains to prove the uniform correction bound. For
\(\abs{v}\geq\rho^{p-1}\), the difference \(\phi_\rho^*(v)-\psi^*(v)\) equals
\(c_p\rho^p\). For \(\abs{v}\leq\rho^{p-1}\), write
\(s:=\abs{v}/\rho^{p-1}\in[0,1]\). Then
\[
    \phi_\rho^*(v)-\psi^*(v)
    =\rho^p\left(\frac12s^2-\frac1q s^q\right).
\]
The function in parentheses has derivative \(s-s^{q-1}\geq0\) on
\([0,1]\), and therefore increases from zero to
\(1/2-1/q=c_p\). This gives~\eqref{eq:bounded-dual-correction}.
In either case, \(\phi_\rho\) is finite, closed, and convex, so
\(\phi_\rho=(\phi_\rho^*)^*\) by biconjugacy, proving
\eqref{eq:scalar-dual-p-huber}.
\end{proof}

\subsubsection{Proof of Proposition~\ref{prop:p-huber-properties}}\label{app:proof-p-huber-properties}

\begin{proof}
By evenness it suffices to consider \(r\geq0\).  At \(r=\rho\), the two
expressions defining \(\phi_\rho\) in \eqref{eq:p-huber-smoothing-surrogate} have the common
value \(\rho^p/2\) and the common derivative
\(\rho^{p-1}\).  Hence \(\phi_\rho\in \mathcal{C}^1(\R)\), with derivative
\(\nabla\phi_\rho(r)=(\max\{\rho,\abs r\})^{p-2}r\).  On
\((0,\rho)\), its second
derivative is \(\rho^{p-2}\); on \((\rho,\infty)\), it is
\((p-1)r^{p-2}\), which belongs to
\([0,(p-1)\rho^{p-2}]\).  Therefore the derivative is nondecreasing and
\(\rho^{p-2}\)-Lipschitz.  The Lipschitz constant is exact because
\(\nabla^2\phi_\rho(r)=\rho^{p-2}\) for \(\abs{r}<\rho\).

Let \(e_\rho(r):=\psi(r)-\phi_\rho(r)\).  For
\(0<r\leq\rho\),
\[
    \nabla e_\rho(r)
    =r^{p-1}-\rho^{p-2}r
    =r\bigl(r^{p-2}-\rho^{p-2}\bigr)\geq0,
\]
because \(p-2<0\), while \(e_\rho(0)=0\).  Thus \(e_\rho\) increases from zero to
\(e_\rho(\rho)=c_p\rho^p\); for \(r>\rho\), it is identically
\(c_p\rho^p\).  Hence
\(0\leq\psi(r)-\phi_\rho(r)\leq c_p\rho^p\) for every
\(r\in\R\), as claimed.

To bound the change in \(\phi_\rho(r)\), fix \(r\) and differentiate with respect to
\(\rho\).  Away from \(\rho=\abs r\),
\[
\partial_\rho\phi_\rho(r)
=
\begin{cases}
-\dfrac{2-p}{2}\rho^{p-1},&0<\rho<\abs r,\\[0.5em]
-\dfrac{2-p}{2}\rho^{p-3}r^2,&\rho>\abs r.
\end{cases}
\]
The two expressions agree at \(\rho=\abs r\), and in either case
\[
    -\frac{2-p}{2}\rho^{p-1}
    \leq \partial_\rho\phi_\rho(r)\leq0.
\]
Therefore, for \(0<\rho_2\leq\rho_1\), integration from
\(\rho_2\) to \(\rho_1\) gives
\[
    0\leq
    \phi_{\rho_2}(r)-\phi_{\rho_1}(r)
    \leq c_p\bigl(\rho_1^p-\rho_2^p\bigr).
\]
Applying Lemma~\ref{lem:coordinatewise-lifting} with
\(L=\rho^{p-2}\) and \(\delta=c_p\rho^p\) establishes the stated
convexity, smoothness, and approximation properties of \(f_\rho\).
Finally, summing the last scalar inequality over the \(m\) residual
coordinates gives the bound in~\eqref{eq:p-huber-radius-drift}.
\end{proof}

\subsubsection{Proof of Proposition~\ref{prop:composite-surrogate-properties}}\label{app:proof-composite-surrogate-properties}

\begin{proof}
By Proposition~\ref{prop:p-huber-properties},
\(f_\rho\) is finite, continuous, and convex.  Since \(\omega\)
is proper, lower semicontinuous, and convex, adding the finite continuous
function \(f_\rho\) shows that \(F_\rho\) is proper,
lower semicontinuous, and convex.

For every \(\vu\in\dom(\omega)\), the regularizer is finite and cancels
from the relevant differences.  Therefore
\[
F(\vu)-F_\rho(\vu)
=f(\vu)-f_\rho(\vu),
\]
so the claimed approximation bound follows from
Proposition~\ref{prop:p-huber-properties}.  Likewise, for
\(0<\rho_2\leq\rho_1\),
\[
F_{\rho_2}(\vu)-F_{\rho_1}(\vu)
=f_{\rho_2}(\vu)-f_{\rho_1}(\vu),
\]
and the bound on \(F_{\rho_2}-F_{\rho_1}\) follows from the same proposition.

Finally, let \(\vu^\star\in\mathcal{U}^\star\).  Since
\(F(\vu^\star)=F^\star\), adding and subtracting
\(f_\rho(\vu)\) and \(f_\rho(\vu^\star)\) gives
\[
\begin{aligned}
F(\vu)-F^\star
&=f(\vu)-f_\rho(\vu)
+F_\rho(\vu)-F_\rho(\vu^\star) \\
&\qquad+f_\rho(\vu^\star)-f(\vu^\star).
\end{aligned}
\]
The claimed inequality follows from
\[
f(\vu)-f_\rho(\vu)\leq c_p m\rho^p
\quad\text{and}\quad
f_\rho(\vu^\star)-f(\vu^\star)\leq0.
\]
\end{proof}

\subsubsection{Proof of Theorem~\ref{thm:optimal-smoothing-surrogate}}

\begin{proof}
Fix \(L>0\) and \(\phi\in\Phi_L\), and let
\(\delta:=\sup_r(\psi(r)-\phi(r))\).  The conclusion is immediate if
\(\delta=+\infty\), so suppose \(\delta<\infty\).  The standard upper model
for a function with \(L\)-Lipschitz derivative gives, for every \(r>0\),
\begin{align*}
\phi(r)
&\leq\phi(0)+\nabla\phi(0)r+\frac{L}{2}r^2,\\
\phi(-r)
&\leq\phi(0)-\nabla\phi(0)r+\frac{L}{2}r^2.
\end{align*}
Since \(\phi(0)\leq\psi(0)=0\), averaging these inequalities and
using the definition of \(\delta\) yields
\[
\delta
\geq
\frac{\psi(r)-\phi(r)+\psi(-r)-\phi(-r)}{2}
\geq \frac{r^p}{p}-\frac{Lr^2}{2}.
\]
The right-hand side is maximized at
\(r=L^{-1/(2-p)}\), where its value is
\(c_pL^{-p/(2-p)}\).  Hence
\(\delta^\star(L)\geq c_pL^{-p/(2-p)}\).

Conversely, choose \(\rho=L^{-1/(2-p)}\).  Then
\(\rho^{p-2}=L\).  The scalar properties established above show that
\(\phi_\rho\) lies below \(\psi\), has \(L\)-Lipschitz derivative,
and has worst-case error
\(c_p\rho^p=c_pL^{-p/(2-p)}\).  Thus
\(\phi_\rho\in\Phi_L\), proving equality and attainment.

For completeness, the same constant has a direct dual derivation on the
convex subclass. If \(\phi\in\Phi_L\) is convex, set \(\mu:=L^{-1}\),
\(g:=\phi^*\), and \(\delta:=\sup_r(\psi(r)-\phi(r))\). Then \(g\) is
\(\mu\)-strongly convex and \(\psi^*\leq g\leq \psi^*+\delta\). At \(p=1\), these inequalities imply \(\dom g=[-1,1]\) and
\[
    \delta\geq\frac{g(1)+g(-1)}{2}
    \geq g(0)+\frac{\mu}{2}
    \geq\frac{\mu}{2}=\frac{1}{2L}.
\]
For \(1<p<2\), strong convexity gives, for every \(a>0\),
\[
    \frac{a^q}{q}+\delta
    \geq \frac{g(a)+g(-a)}{2}
    \geq g(0)+\frac{\mu a^2}{2}
    \geq \frac{\mu a^2}{2}.
\]
Therefore
\[
    \delta\geq \frac{\mu a^2}{2}-\frac{a^q}{q}.
\]
Maximizing the right-hand side at \(a=\mu^{1/(q-2)}\) gives
\[
    \delta\geq
    \left(\frac12-\frac1q\right)\mu^{q/(q-2)}
    =c_pL^{-p/(2-p)}.
\]
For \(\mu=\rho^{2-p}\), the maximizing dual point is
\(a=\rho^{p-1}\), exactly the breakpoint in
\eqref{eq:p-huber-conjugate}.
\end{proof}

\subsection{Proofs for Section~\ref{sec:optimal-matrix-vector-complexity}}

\subsubsection{Proof of Lemma~\ref{lem:accelerated-drift}}
\label{app:proof-accelerated-drift}

\begin{proof}
The proximal subproblem is strongly convex, so its minimizer is unique and
belongs to \(\dom(\omega)\).  Since \(\dom(\omega)\) is convex and
\(\vy_0=\vz_0=\vu_0\in\dom(\omega)\), all iterates
\(\vx_t,\vy_t,\vz_t\) lie in \(\dom(\omega)\).

Fix \(t\) and \(\vu\in\dom(\omega)\), and write
\(\vg_t:=\nabla f_t(\vx_t)\).  Smoothness of \(f_t\) and
\(\vy_t-\vx_t=\alpha_t(\vz_t-\vz_{t-1})\) give
\begin{align}
f_t(\vy_t)
&\leq
f_t(\vx_t)
+\alpha_t\ip{\vg_t}{\vz_t-\vz_{t-1}}
+\frac{L_t\alpha_t^2}{2}\norm{\vz_t-\vz_{t-1}}^2.
\label{eq:app-composite-smooth-model}
\end{align}
Convexity of \(f_t\), applied at \(\vx_t\) to
\(\vy_{t-1}\) and \(\vu\), yields
\begin{equation}
f_t(\vx_t)
+\alpha_t\ip{\vg_t}{\vu-\vz_{t-1}}
\leq
(1-\alpha_t)f_t(\vy_{t-1})
+\alpha_t f_t(\vu).
\label{eq:app-composite-convex-loss}
\end{equation}
Convexity of \(\omega\) also gives
\begin{equation}
    \omega(\vy_t)
    \leq(1-\alpha_t)\omega(\vy_{t-1})
       +\alpha_t\omega(\vz_t).
\label{eq:app-composite-convex-omega}
\end{equation}

Proximal optimality provides \(\vs_t\in\partial\omega(\vz_t)\) such that
\[
    \vg_t+\vs_t+\gamma_t(\vz_t-\vz_{t-1})=\vzero.
\]
The subgradient inequality and the polarization identity therefore imply
\begin{align}
&\ip{\vg_t}{\vz_t-\vu}+\omega(\vz_t)-\omega(\vu)
\nonumber\\
&\qquad\leq
\frac{\gamma_t}{2}
\left(
\norm{\vu-\vz_{t-1}}^2
-\norm{\vu-\vz_t}^2
-\norm{\vz_t-\vz_{t-1}}^2
\right).
\label{eq:app-proximal-three-point}
\end{align}
Combining
\eqref{eq:app-composite-smooth-model}--
\eqref{eq:app-proximal-three-point} gives
\begin{align}
F_t(\vy_t)-F_t(\vu)
&\leq
(1-\alpha_t)
\bigl(F_t(\vy_{t-1})-F_t(\vu)\bigr)
\nonumber\\
&\quad+
\frac{\alpha_t\gamma_t}{2}
\left(
\norm{\vu-\vz_{t-1}}^2-\norm{\vu-\vz_t}^2
\right)
\nonumber\\
&\quad-
\frac{\alpha_t}{2}(\gamma_t-\alpha_tL_t)
\norm{\vz_t-\vz_{t-1}}^2.
\label{eq:app-composite-one-step}
\end{align}
The last term is nonpositive. Multiplying \eqref{eq:app-composite-one-step} by \(\Gamma_t\) and using
\(\Gamma_{t-1}=(1-\alpha_t)\Gamma_t\) and
\(\beta_t=\alpha_t\gamma_t\Gamma_t\) yields
\begin{align}
\Gamma_t\bigl(F_t(\vy_t)-F_t(\vu)\bigr)
&\leq
\Gamma_{t-1}\bigl(F_t(\vy_{t-1})-F_t(\vu)\bigr)
\nonumber\\
&\quad+
\frac{\beta_t}{2}
\left(
\norm{\vu-\vz_{t-1}}^2-\norm{\vu-\vz_t}^2
\right).
\label{eq:app-composite-same-surrogate}
\end{align}
Finally, \(\rho_t\leq\rho_{t-1}\) and the bound on \(F_t-F_{t-1}\)
in Proposition~\ref{prop:composite-surrogate-properties} give
\[
F_t(\vy_{t-1})-F_t(\vu)
\leq
F_{t-1}(\vy_{t-1})-F_{t-1}(\vu)
+c_p m(\rho_{t-1}^p-\rho_t^p).
\]
Substituting this into \eqref{eq:app-composite-same-surrogate} proves
\eqref{eq:accelerated-drift}.
\end{proof}

\subsubsection{Proof of Theorem~\ref{thm:upper-bound}}
\label{app:proof-fixed-radius}

\begin{proof}
Choose \(\vu^\star\in\mathcal{U}^\star\) with
\(\norm{\vu_0-\vu^\star}\leq D\).  For the fixed radius, the drift term in
Lemma~\ref{lem:accelerated-drift} vanishes.  Moreover,
\[
\gamma_t=\frac{2L_T}{t}
\geq\frac{2L_T}{t+1}=\alpha_tL_T,
\qquad
\beta_t
=\frac{2}{t+1}\frac{2L_T}{t}\frac{t(t+1)}{2}
=2L_T.
\]
Summing \eqref{eq:accelerated-drift} from \(t=1\) to \(T\), using
\(\Gamma_0=0\), \(\vz_0=\vu_0\), and the constancy of \(\beta_t\), gives
\[
\Gamma_T\bigl(
F_{\rho_T}(\vy_T)-
F_{\rho_T}(\vu^\star)
\bigr)
\leq
L_T\bigl(
\norm{\vu^\star-\vu_0}^2-
\norm{\vu^\star-\vz_T}^2
\bigr)
\leq L_TD^2.
\]
Combining this estimate with the objective-gap bound in
Proposition~\ref{prop:composite-surrogate-properties} and
\(\Gamma_T=T(T+1)/2\) yields, for every fixed \(\rho>0\),
\begin{equation*}
F(\vy_T)-F^\star
\leq
c_p m\rho^p
+\frac{2\norm{\mA}^2D^2}{T(T+1)}\rho^{p-2}.
\end{equation*}
The derivative of the right-hand side vanishes uniquely at
\[
    \rho^2=\frac{4\norm{\mA}^2D^2}{mT(T+1)},
\]
which is \eqref{eq:fixed-radius}.  Substituting this radius into the two first-order balance relations gives
\[
F(\vy_T)-F^\star
\leq\frac{m}{p}\rho_T^p
=\frac{2^p}{p}m^{\frac{2-p}{2}}
 \left(\frac{\norm{\mA}D}{\sqrt{T(T+1)}}\right)^p.
\]
The final inequality in \eqref{eq:fixed-rate} follows from
\(T(T+1)\geq T^2\).  Each iteration evaluates
\(\mA\vx_t\), one adjoint product, and one proximal map, establishing the
stated counts.
\end{proof}

\subsubsection{Proof of Theorem~\ref{thm:variable-radius}}
\label{app:proof-variable-radius}

\begin{proof}
Let \(a_t:=\norm{\vz_t-\vu^\star}^2\).  Under
\(\rho_t=\bar\rho/(t+1)\) and
\(\gamma_t=\alpha_tL_t\),
\begin{equation*}
\beta_t
=
\alpha_t^2L_t\Gamma_t
=
\frac{2t}{t+1}\norm{\mA}^2\bar\rho^{p-2}(t+1)^{2-p}
=
2\norm{\mA}^2\bar\rho^{p-2}t(t+1)^{1-p}.
\end{equation*}
The function \(s\mapsto s(s+1)^{1-p}\) is increasing on
\((0,\infty)\), because its derivative is
\((s+1)^{-p}[1+(2-p)s]>0\).  Hence \(\{\beta_t\}\) is non-decreasing.

Summing Lemma~\ref{lem:accelerated-drift} with
\(\vu=\vu^\star\) gives
\begin{align*}
\Gamma_T
\bigl(F_T(\vy_T)-F_T(\vu^\star)\bigr)
&\leq
\frac12\sum_{t=1}^T\beta_t(a_{t-1}-a_t)
+c_p m\sum_{t=1}^T
\Gamma_{t-1}(\rho_{t-1}^p-\rho_t^p).
\end{align*}
By Abel summation and \(a_t\leq R^2\),
\begin{align*}
\sum_{t=1}^T\beta_t(a_{t-1}-a_t)
&=
\beta_1a_0
+\sum_{t=1}^{T-1}(\beta_{t+1}-\beta_t)a_t
-\beta_Ta_T
\leq
\beta_TR^2.
\end{align*}
The uniform approximation bound therefore yields
\begin{align}
F(\vy_T)-F^\star
&\leq
\frac{\beta_TR^2}{2\Gamma_T}
+c_p m\rho_T^p
+\frac{c_p m}{\Gamma_T}
\sum_{t=1}^T
\Gamma_{t-1}(\rho_{t-1}^p-\rho_t^p).
\label{eq:app-variable-prerate}
\end{align}
The first term is
\begin{equation*}
\frac{\beta_TR^2}{2\Gamma_T}
=
\frac{2\norm{\mA}^2R^2\bar\rho^{p-2}}{(T+1)^p},
\end{equation*}
and the second is
\(c_p m\bar\rho^p(T+1)^{-p}\).
For the drift term, the mean-value theorem gives
\begin{equation*}
 t^{-p}-(t+1)^{-p}\leq p t^{-p-1},
\end{equation*}
and \(\Gamma_{t-1}\leq t^2/2\).  Since
\(1\leq p<2\),
\begin{equation*}
\sum_{t=1}^T t^{1-p}
\leq
\frac{(T+1)^{2-p}}{2-p}.
\end{equation*}
Consequently,
\begin{align*}
\frac{c_p m}{\Gamma_T}
\sum_{t=1}^T
\Gamma_{t-1}(\rho_{t-1}^p-\rho_t^p)
&\leq
\frac{c_p mp\bar\rho^p}{2\Gamma_T}
\sum_{t=1}^Tt^{1-p}
\nonumber\\
&\leq
\frac{m\bar\rho^p}{2}
\frac{(T+1)^{1-p}}{T}
\leq
m\bar\rho^p(T+1)^{-p}.
\end{align*}
Substitution into \eqref{eq:app-variable-prerate} gives the two-term bound
\begin{equation}
F(\vy_T)-F^\star
\leq
\frac{1}{(T+1)^p}
\left[
2\norm{\mA}^2R^2\bar\rho^{p-2}
+\frac{p+2}{2p}m\bar\rho^p
\right].
\label{eq:variable-two-term}
\end{equation}
It remains to optimize over \(\bar\rho\).  Differentiating the bracket in
\eqref{eq:variable-two-term} gives
\begin{equation*}
2\norm{\mA}^2R^2(p-2)\bar\rho^{p-3}
+\frac{p+2}{2}m\bar\rho^{p-1}=0,
\end{equation*}
whose unique positive solution is
\begin{equation*}
\bar\rho^2
=
\frac{4(2-p)}{(p+2)m}\norm{\mA}^2R^2.
\end{equation*}
At this value, the bracket equals
\begin{equation*}
\frac{4}{p}
\left(\frac{p+2}{4(2-p)}\right)^{\frac{2-p}{2}}
 m^{\frac{2-p}{2}}(\norm{\mA}R)^p,
\end{equation*}
which proves \eqref{eq:variable-rate}.  The oracle count is the same as in the
fixed-smoothing variant.
\end{proof}

\subsubsection{Proof of Theorem~\ref{thm:matvec-lower-bound}}
\label{app:proof-matvec-lower-bound}

The proof has two ingredients.  We adapt the rotating hard-instance construction of \citet{ouyang2021lower}, which bounds the Euclidean residual of an arbitrary deterministic method (Lemma~\ref{lem:fixed-rhs-euclidean-hard}).

\begin{lemma}
\label{lem:fixed-rhs-euclidean-hard}
Let \(Q \geq 1\) be an integer, set \(\ell=2Q+4\) and \(d = 2\ell+1\), and define
\begin{equation*}
    \mB=
    \begin{pmatrix}
        1 &        &        &   \\
       -1 & 1      &        &   \\
          & \ddots & \ddots &   \\
          &        & -1     & 1
    \end{pmatrix}\in\R^{2\ell\times2\ell},
    \qquad
    \mC=\begin{pmatrix}\mB&\vzero\\\vzero^\top&2\end{pmatrix},
    \qquad
    \vc=\begin{pmatrix}\vone_{2\ell}\\0\end{pmatrix}.
\end{equation*}
For every deterministic method given \(\vc,Q\) and accessing
an unknown matrix only through at most \(Q\) adaptive queries
to the paired oracle
\begin{equation}\label{eq:pairedoracle}
(\vxi,\vnu)\mapsto(\widetilde{\mA}\vxi,\widetilde{\mA}^{\top}\vnu),
    \qquad \vxi,\vnu\in\R^d,
\end{equation}
there exist orthogonal \(\mU,\mV \in \R^{d\times d}\) with \(\mV\vc = \vc\) such that, on the fixed matrix \(\widetilde{\mA} = \mV^\top \mC \mU\), the method's output $\vxi_Q$ satisfies
\begin{equation}
    \norm{\widetilde{\mA} \vxi_Q - \vc} \geq \sqrt{\ell}.
\end{equation}
\end{lemma}

\begin{proof}
Throughout the proof, the subscript \((_\circ)\) distinguishes the notation of Ouyang and Xu from ours.

Fix a deterministic method \(\mathcal N\) in the oracle model~\eqref{eq:pairedoracle}; if it terminates before \(Q\), the remaining queries can be padded with zero without changing its output.
The oracle in \citet[Equation (1.7)]{ouyang2021lower} returns \((\nabla f_{\circ}, \mA\vxi, \mA^\top \vv)\), and thus our paired oracle~\eqref{eq:pairedoracle} is a specialization of their oracle obtained by setting \(f_\circ\equiv 0\) and discarding the gradient component \(\nabla f_{\circ}\).

In their construction (2.2)--(2.3), take \(m_{\circ}=n_{\circ}=d\), \(k_{\circ}=\ell\), \(\mG_{\circ} = 2\). Choosing \(L_{A\circ}=2\) in their (2.20) gives our right-hand side \(\vc\). Note that our matrix \(\mB\) is obtained by reversing the row order of their \(\mB_\circ\), which does not affect the analysis.

We can therefore apply \citet[Proposition 3.1]{ouyang2021lower} with \(t_\circ = Q\), \(\mH = \vzero\), and \(\vh = \vzero\), taking both domains to be \(\R^d\). Its conditions hold since \(\ell < d/2\), \(Q\leq \ell/2-1\), and \(\mH_\circ\) vanishes. The proposition yields orthogonal \(\mU,\mV\) with \(\mV \vc = \vc\) and \(\widetilde{\mA} = \mV^\top \mC \mU\) such that
\[
    \norm{\widetilde{\mA} \vxi_Q - \vc} \geq \min_{{\vw} \in \mathcal{K}_{\ell-1}} \norm{\mC\vw -\vc},
\]
where \(\vxi_Q\) is the method's output on this fixed matrix and \(\mathcal{K}_{\ell-1}\) is the Krylov subspace that, by \citet[Lemma 2.2]{ouyang2021lower}, satisfies
\begin{equation*}
    \mathcal{K}_{\ell-1} = \{
    \vw \in \R^d: w_1 = \cdots = w_\ell = w_d = 0
    \}.
\end{equation*}
For every \(\vw \in \mathcal{K}_{\ell-1}\), the first \(\ell\) entries of \(\mC\vw\) are therefore zero, while the corresponding entries of \(\vc\) are one. Hence \(\norm{\mC\vw - \vc}^2 \geq \ell\), which proves the claim.
\end{proof}

The previous lemma gives a Euclidean residual bound in \(\R^d\).
We now transfer it to the \(\ell_p\) norm in \(\R^m\), using the
conclusion following the proof of Theorem~4 in
\citet{kashin1977diameters}.

\begin{lemma}\label{lem:euclidean-section-embedding}
There exists a universal constant \(c_0\in(0,1)\) such that, for all integers \(m,d\geq1\) with \(d\leq m/2\), there is a matrix \(\mR\in\R^{m\times d}\) satisfying \(\mR^{\top}\mR=\mI_d\) and, simultaneously for every \(\vs\in\R^d\) and every \(p\in[1,2]\),
\begin{equation}
\label{eq:lower-R-lp}
    \norm{\mR\vs}_p^p
    \geq
    {c_0}^p m^{1-p/2}\norm{\vs}^p.
\end{equation}
\end{lemma}

\begin{proof}

Similarly, a subscript or superscript \((\dagger)\) distinguishes Kashin's notation from ours. The result of \citet{kashin1977diameters} states that, for every \(\theta_\dagger > 0\), there exists a constant \(C_{\theta}^\dagger > 0\) such that, for every integer \(n_\dagger \geq 1\), there exists a subspace \(L_{n_\dagger, \theta_\dagger} \subseteq \R^{n_\dagger}\) with \(\dim{L_{n_\dagger, \theta_\dagger}} \geq n_\dagger (1-\theta_\dagger)\). Every \(\rvx \in L_{n_\dagger, \theta_\dagger}\) satisfies \(C_{\theta}^\dagger \norm{\rvx}_2 \leq \frac{1}{\sqrt{n_\dagger}} \norm{\rvx}_1 \leq \norm{\rvx}_2\).

To apply the result, set \(n_\dagger = m\) and \(\theta_\dagger = 1/2\). Then there exists a subspace \(\mathcal{E}\subseteq \R^m\) such that \(\dim \mathcal{E} \geq  m/2\) and for all \(\vv \in \mathcal{E}\),
\begin{equation}\label{eq:v1norm}
    \norm{\vv}_1 \geq C_{1/2}^{\dagger} \sqrt{m} \norm{\vv}_2.
\end{equation}
Since \(d\leq m/2\leq\dim\mathcal E\), choose \(d\) orthonormal vectors in \(\mathcal E\) as the columns of \(\mR\). Then \(\mR^\top\mR=\mI_d\) and
\(\norm{\mR\vs}_2=\norm{\vs}_2\).
Since \(\mR\vs \in \mathcal{E}\), \eqref{eq:v1norm} gives
\begin{equation}
    \norm{\mR\vs}_1 \geq C_{1/2}^{\dagger} \sqrt{m} \norm{\vs}_2.
\end{equation}

We now extend the bound to the \(\ell_p\) norm for all \(p \in [1,2]\). By H\"older's inequality \(\norm{\vv}_p \geq m^{1/p -1} \norm{\vv}_1\),
substituting \(\vv = \mR\vs\) gives
\begin{equation}
    \norm{\mR\vs}_p \geq m^{1/p -1} \norm{\mR\vs}_1 \geq C_{1/2}^{\dagger} ~m^{1/p -1/2}\norm{\vs}_2 .
\end{equation}
Raising both sides to the power \(p\) and letting \(c_0 = \min\{C_{1/2}^{\dagger},1/2\}\) proves the claim.
\end{proof}

We now combine the two lemmas to prove the lower bound.
\begin{proof}[Proof of Theorem~\ref{thm:matvec-lower-bound}]
Let \(\ell=2Q+4\), \(d=2\ell+1\), and take \(\mC,\vc\) from
Lemma~\ref{lem:fixed-rhs-euclidean-hard}.
Since \(m\geq2d\), choose \(\mR\) from
Lemma~\ref{lem:euclidean-section-embedding}. Since \(n\geq d+1\),
choose \(\mE\in\R^{n\times d}\) satisfying
\[
    \mE^\top\mE=\mI_d,
    \qquad \mE^\top\vu_0=\vzero.
\]
The solution of \(\mC\vxi=\vc\) is
\[
    \vxi_0^\star=(1,2,\ldots,2\ell,0)^\top,
    \qquad
    \norm{\vxi_0^\star}^2
    =\frac{\ell(2\ell+1)(4\ell+1)}{3}.
\]
Set
\[
    a:=\frac{D}{\norm{\vxi_0^\star}},
    \qquad
    \vb:=\frac{a\Lambda}{2}\mR\vc,
    \qquad \omega\equiv0.
\]
These quantities are fixed independently of the method and the hard matrix.

Fix any deterministic method \(\mathcal M\) in the theorem.
For \(\mA=(\Lambda/2)\mR\widetilde{\mA}\mE^\top\), the identity
\[
    (\mA\vu,\mA^\top\vv)
    =\frac{\Lambda}{2}
      \bigl(\mR\widetilde{\mA}(\mE^\top\vu),
            \mE\widetilde{\mA}^\top(\mR^\top\vv)\bigr)
\]
shows that each paired query to \(\mA\) requires only one paired
query to \(\widetilde{\mA}\).
Since the other transformations are fixed and the proximal map is the
identity, \(\vxi_Q:=a^{-1}\mE^\top\vu_Q\) is the output of a
deterministic method using at most \(Q\) queries to
\(\widetilde{\mA}\). Apply
Lemma~\ref{lem:fixed-rhs-euclidean-hard} to this induced method
and choose the corresponding
\[
    \widetilde{\mA}=\mV^\top\mC\mU,
    \qquad
    \mA=\frac{\Lambda}{2}\mR\widetilde{\mA}\mE^\top.
\]
Then \(\mV\vc=\vc\) and
\(\norm{\widetilde{\mA}\vxi_Q-\vc}\geq\sqrt\ell\), where
\(\vu_Q\) is the output of \(\mathcal M\) on this instance.

For every \(\vw\in\R^{2\ell}\),
\[
    \norm{\mB\vw}^2
    =w_1^2+\sum_{i=2}^{2\ell}(w_i-w_{i-1})^2
    \leq4\norm{\vw}^2,
\]
so \(\norm{\mC}=\max\{\norm{\mB},2\}=2\).
The orthonormal columns of \(\mR,\mE\) and the identity
\(\mR^\top\mA\mE=(\Lambda/2)\widetilde{\mA}\) give
\[
    \norm{\mA}
    =\frac{\Lambda}{2}\norm{\widetilde{\mA}}
    =\frac{\Lambda}{2}\norm{\mC}=\Lambda.
\]
Since \(\mC\) is invertible and \(\mV\vc=\vc\),
\(\widetilde{\mA}^{-1}\vc=\mU^\top\vxi_0^\star\).
Thus \(\vu^\star:=\vu_0+a\mE\mU^\top\vxi_0^\star\)
satisfies \(\mA\vu^\star=\vb\) and
\(\norm{\vu^\star-\vu_0}=D\), so \(F^\star=0\).
Every minimizer \(\vu\) satisfies
\(\mE^\top\vu=a\mU^\top\vxi_0^\star\), and hence
\[
    \norm{\vu-\vu_0}
    \geq\norm{\mE^\top(\vu-\vu_0)}
    =a\norm{\vxi_0^\star}=D.
\]
Equality holds at \(\vu^\star\), proving the required initial distance.

Finally,
\[
    \mA\vu_Q-\vb
    =\frac{a\Lambda}{2}\mR
      (\widetilde{\mA}\vxi_Q-\vc).
\]
The two lemmas therefore yield
\begin{align*}
    F(\vu_Q)-F^\star
    &=\frac1p\left(\frac{a\Lambda}{2}\right)^p
        \norm{\mR(\widetilde{\mA}\vxi_Q-\vc)}_p^p\\
    &\geq\frac{c_0^p}{p}m^{1-p/2}
        \left(\frac{a\Lambda\sqrt\ell}{2}\right)^p=\frac{c_0^p}{p}m^{1-p/2}
        \left(\frac{\sqrt3\,\Lambda D}
        {2\sqrt{(2\ell+1)(4\ell+1)}}\right)^p\\
    &\geq\frac1p\left(\frac{\sqrt3\,c_0}{4\sqrt2}\right)^p
        m^{1-p/2}\left(\frac{\Lambda D}{2Q+5}\right)^p,
\end{align*}
where the last step uses
\((2\ell+1)(4\ell+1)\leq8(\ell+1)^2\) and \(\ell+1=2Q+5\).
\end{proof}

\section{Additional discussion}

\subsection{Fixed smoothing via restart}
\label{app:variants}

The fixed-smoothing guarantee in Theorem~\ref{thm:upper-bound} uses a radius tuned
to the terminal horizon.  When the horizon is not known in advance, a
restart-doubling scheme removes this requirement without changing the order
of the rate.  This construction runs the fixed-smoothing variant with horizon-tuned radii
at geometrically increasing stage lengths.  It is distinct from the
variable-smoothing variant in Theorem~\ref{thm:variable-radius},
which uses a single run with \(\rho_t=\bar\rho/(t+1)\).

\begin{proposition}
\label{prop:fixed-doubling}
For stages \(k=0,1,2,\ldots\), set \(T_k:=2^k\) and run the fixed-smoothing
variant of Theorem~\ref{thm:upper-bound} for \(T_k\) iterations from the common
initial point \(\vu_0\), using the radius tuned to \(T_k\).  After a total iteration
budget \(N\geq1\), denote the output of the last completed stage by \(\vu_N\)
and return it.  Then
\begin{equation*}
    F(\vu_N)-F^\star
    \leq
    \frac{8^p}{p}\,
    m^{\frac{2-p}{2}}
    \left(\frac{\norm{\mA}D}{N}\right)^p.
\end{equation*}
The resulting scheme requires neither a prescribed terminal horizon nor any bounded-iterate
assumption.
\end{proposition}

\begin{proof}
Let \(K\) be the last completed stage.  Since
\(\sum_{k=0}^{K}2^k=2^{K+1}-1\leq N\), while stage \(K+1\) is not completed,
we have \(N<2^{K+2}-1\), and therefore \(T_K=2^K\geq N/4\).  Applying
Theorem~\ref{thm:upper-bound} to stage \(K\) gives
\begin{equation*}
F(\vu_N)-F^\star
\leq
\frac{2^p}{p}m^{\frac{2-p}{2}}
\left(\frac{\norm{\mA}D}{T_K}\right)^p
\leq
\frac{8^p}{p}m^{\frac{2-p}{2}}
\left(\frac{\norm{\mA}D}{N}\right)^p.
\end{equation*}
\end{proof}

\subsection{Constrained variable smoothing}
\label{app:projected-variable-radius}

For a general nonsmooth \(\omega\), the bounded-domain variant uses the
proximal map of the sum of \(\omega\) and the ball indicator.

\begin{corollary}[Variable smoothing without a bounded-iterate assumption]
\label{cor:projected-variable-radius}
Let
\(\mathcal B_D:=\{\vu:\norm{\vu-\vu_0}\leq D\}\), and assume that the
constrained proximal map of
\(\omega+\iota_{\mathcal B_D}\) is available.  Replace the \(\vz\)-update
in Algorithm~\ref{alg:hs-apg} by
\begin{equation}
\vz_t
:=\prox_{(\omega+\iota_{\mathcal B_D})/\gamma_t}
\left(\vz_{t-1}-\gamma_t^{-1}\vg_t\right).
\label{eq:constrained-prox-update}
\end{equation}
Run the variable schedule of Theorem~\ref{thm:variable-radius} with
\(R=2D\).  Then
\begin{equation}
\label{eq:projected-variable-rate}
    F(\vy_T)-F^\star
    \leq
    2^pC_p\,
    m^{\frac{2-p}{2}}
    \left(\frac{\norm{\mA}D}{T+1}\right)^p.
\end{equation}
The method uses exactly \(2T\) individual matrix--vector products and
\(T\) evaluations of the constrained proximal map.
\end{corollary}

\begin{proof}
Choose a nearest minimizer \(\vu^\star\in\mathcal{U}^\star\).  Then
\(\norm{\vu^\star-\vu_0}\leq D\), so \(\vu^\star\in\mathcal B_D\).  The
constrained problem
\[
\min_{\vu}\ f(\vu)+\omega(\vu)+\iota_{\mathcal B_D}(\vu)
\]
has the same optimal value \(F^\star\), because its feasible set contains the global
minimizer \(\vu^\star\).  Update \eqref{eq:constrained-prox-update} is
exactly Algorithm~\ref{alg:hs-apg} applied to the regularizer
\(\omega+\iota_{\mathcal B_D}\), so the composite three-point inequality
in the proof of Lemma~\ref{lem:accelerated-drift} remains valid.

Every \(\vz_t\) lies in \(\mathcal B_D\), and every \(\vy_t\), being a
convex combination of points in this ball, also lies there.  Hence
\[
    \norm{\vz_t-\vu^\star}
    \leq\norm{\vz_t-\vu_0}+\norm{\vu_0-\vu^\star}
    \leq2D.
\]
Theorem~\ref{thm:variable-radius} applies with \(R=2D\), and its objective
gap is exactly \(F(\vy_T)-F^\star\), because
\(\vy_T\in\mathcal B_D\).  This proves
\eqref{eq:projected-variable-rate}; the constrained proximal evaluation
requires no matrix--vector product.
\end{proof}

\subsection{Backtracking line search}\label{app:backtracking-line-search}

\begin{proposition}\label{prop:backtracking-variable}
Consider the radius schedule and bounded-sequence assumption of
Theorem~\ref{thm:variable-radius}.  At iteration \(t\), accept a local
estimate \(\widehat L_t\) for which the test
\eqref{eq:p-huber-line-search-test} holds, set
\(\gamma_t=\alpha_t\widehat L_t\), and impose
\begin{equation}
0<\widehat L_1\leq L_1,
\qquad
\frac{t^2-1}{t^2}\widehat L_{t-1}
\leq\widehat L_t\leq L_t
\quad(t\geq2).
\label{eq:line-search-beta-floor}
\end{equation}
Then the conclusion \eqref{eq:variable-rate} remains valid.  A backtracking
search capped at \(L_t\) always terminates.
\end{proposition}

\begin{proof}
The accepted test is precisely the smooth upper-model inequality needed in
\eqref{eq:app-composite-smooth-model}, evaluated at the actual pair
\((\vx_t,\vy_t)\).  All remaining steps in the proof of
Lemma~\ref{lem:accelerated-drift} use only convexity and proximal optimality,
so the lemma remains valid with
\[
    \widehat\beta_t
    :=\alpha_t^2\widehat L_t\Gamma_t
    =\frac{2t}{t+1}\widehat L_t.
\]
The lower safeguard in \eqref{eq:line-search-beta-floor} is equivalent to
\(\widehat\beta_t\geq\widehat\beta_{t-1}\).  Consequently, the Abel
summation argument in the proof of Theorem~\ref{thm:variable-radius} applies
unchanged.  Since \(\widehat L_T\leq L_T\), we also have
\(\widehat\beta_T\leq\beta_T\), so the same worst-case bound follows.

The global constant \(L_t\) satisfies the test by
Proposition~\ref{prop:p-huber-properties}.  Moreover, for the
variable schedule, \(L_t>L_{t-1}\), while the lower safeguard is strictly
smaller than \(\widehat L_{t-1}\leq L_{t-1}\).  Thus the admissible interval
in \eqref{eq:line-search-beta-floor} contains \(L_t\), and a search capped at
that value terminates.
\end{proof}

Each line-search trial evaluates one candidate proximal point and
\(f_t(\vy_t)\); the latter requires one additional forward product.
Rejected trials therefore add one constrained or unconstrained proximal
evaluation and one forward product each, but no additional adjoint product.
Without the lower safeguard in \eqref{eq:line-search-beta-floor}, every
accepted step still has a certified local upper model, but the global Abel
summation argument need not hold.

\subsection{Least-squares initialization and residual-based scale selection}
\label{app:ols-variant}

For \(\omega\equiv0\), we compute a least-squares solution
and residual
\begin{equation*}
\vu_{\mathrm{OLS}}\in\arg\min_{\vu}\tfrac12\|\mA\vu-\vb\|_2^2,
\qquad
\vr_{\mathrm{OLS}}=\mA\vu_{\mathrm{OLS}}-\vb.
\end{equation*}
For nonzero residual, we initialize Algorithm~\ref{alg:hs-apg} at
\(\vu_{\mathrm{OLS}}\) and use
\begin{equation*}
\bar\rho_{\mathrm{OLS}}=2\|\vr_{\mathrm{OLS}}\|_\infty,
\qquad
\rho_t=\frac{\bar\rho_{\mathrm{OLS}}}{t+1},
\qquad
L_t=\|\mA\|^2\rho_t^{p-2},
\qquad
\gamma_t=\alpha_tL_t.
\end{equation*}
There is no backtracking. In exact arithmetic, a zero residual already
solves the unregularized problem.

The factor two makes \(\rho_1=\|\vr_{\mathrm{OLS}}\|_\infty\), so all
initial residual coordinates lie in the quadratic part of the first
surrogate. For an exact least-squares solution, the normal equations give
\begin{equation*}
\nabla f_{\rho_1}(\vu_{\mathrm{OLS}})
=\rho_1^{p-2}\mA^\top\vr_{\mathrm{OLS}}=\vzero.
\end{equation*}
Convexity therefore makes \(\vu_{\mathrm{OLS}}\) a minimizer of the first
surrogate. The algorithm then decreases the smoothing radius in subsequent
iterations \citep{bubeck2018homotopy}.

The proof of Theorem~\ref{thm:variable-radius} gives the two-term estimate
\eqref{eq:variable-two-term} for every \(\bar\rho>0\). Setting
\(\bar\rho=\bar\rho_{\mathrm{OLS}}\) gives an
\(\mathcal O(T^{-p})\) objective-gap rate if the accelerated iterates satisfy
\(\|\vz_t-\vu^\star\|\leq R\) with a common \(R\) independent of \(T\).
The radius is determined by the least-squares residual, so running the
method requires neither \(R\) nor a reference solution. An approximate
least-squares solution with positive residual scale retains this conditional
rate; exact least squares is required for the stationarity identity above.

\section{Experimental details}\label{app:experiment-details}

\subsection{Algorithm settings and evaluation}
\label{app:method-details}

\paragraph{Smoothing methods.}
The fixed-smoothing, variable-smoothing, and backtracking methods apply
Algorithm~\ref{alg:hs-apg}. In the two unregularized studies, they start from
\(\vu_0=\mathbf 0\) with \(\omega\equiv0\). In the regularized and constrained
studies, they retain \(\omega\) through its proximal operator. Fixed smoothing
uses the radius in \eqref{eq:fixed-radius}, and variable smoothing uses the
schedule in \eqref{eq:variable-rbar}. In the controlled study, fixed smoothing
uses \(T=30{,}000\) and the known initial distance \(D\), while variable
smoothing and backtracking use \(R=D\) as a distance proxy. The initial
points, horizons, and distance choices for the other
studies are specified in their respective setups. The least-squares
initialization variant used in the real-data study is described in
Appendix~\ref{app:realdata_distance_free}.

Backtracking uses the variable-smoothing schedule and tests
\eqref{eq:p-huber-line-search-test} with a local curvature estimate. The first
trial is \(L_1/2\); for \(t\geq2\), it starts from
\[
\min\left\{L_t,\frac12\widehat L_{t-1}
\left(\frac{\rho_t}{\rho_{t-1}}\right)^{p-2}\right\}.
\]
A failed trial doubles the estimate, capped at \(L_t\). This implementation
omits the lower bound in \eqref{eq:line-search-beta-floor} required by
Proposition~\ref{prop:backtracking-variable}.

\paragraph{Baseline methods.}
UFGM follows the Euclidean method of \citet{nesterov2015universal}. The
unregularized studies use \(\omega\equiv0\), while the regularized and
constrained studies include the proximal operator of \(\omega\). The initial
curvature trial is \(L_0=1\), and failed model tests double the trial estimate.

GAGD+ follows \citet{diakonikolas2024complementary}. The unregularized studies
use the non-composite \(\ell_p\)-regression specialization with its native
\(\ell_p\) geometry and no supplied induced operator norm. The regularized
and constrained studies use the Euclidean specialization with zero
uniform-convexity parameter and cumulative regularization weights. The
initial curvature trial is \(M_0=1\), and failed model tests double the trial
estimate.

For UFGM and GAGD+, the prescribed accuracy enters the recurrence and is set
to \(\epsilon=\tau(F(\vu_0)-F_{\mathrm{ref}})\) before each run. The
controlled study uses \(\tau=10^{-9}\), with both reported thresholds
evaluated along that run. The real-data, regularized, and constrained studies
use separate runs for \(\tau=10^{-4}\) and \(\tau=10^{-6}\), and each target
is evaluated using its matching configuration.

AGDA uses the unconstrained, target-free method of \citet{ren2026accelerated}
in the two unregularized studies, with \(\bar r_0=10^{-2}\) and
\(\beta_0=10^{-3}\), and returns the iterate with the lowest objective value
so far. With iterations indexed from \(k=0\), the implementation uses the
bisection tolerance \(\beta_0/[2(k+1)^2]\) in place of \(\beta_0/(2k^2)\) in
their Algorithm~5. Bisection retains the upper endpoint satisfying the
line-search condition and uses a machine-precision lower bound on the
interval-width tolerance. Baseline initialization parameters are held fixed
across instances.

\paragraph{Computational cost and evaluation.}
Each multiplication by \(\mA\) or \(\mA^\top\) counts as one product. The
count includes initialization and all line-search trials, including trials
interrupted by the budget limit. Products used solely for evaluation and
reporting are recorded separately and excluded. A previously computed
product can be reused exactly without increasing the count. For fixed and
variable smoothing without backtracking, each completed iteration uses one
forward and one adjoint product; a completed backtracking iteration with
\(N_t^{\mathrm{rej}}\) rejected trials uses \(3+N_t^{\mathrm{rej}}\) products.

Evaluation uses the original nonsmoothed objective at each method's
prescribed output after every completed iteration. The cost of reaching a
target is the cumulative product count at the first output satisfying it.
The main comparisons in the two unregularized studies use a budget of \(60{,}000\) products; Appendix~\ref{app:realdata-budget} examines a larger budget on selected real-data problems. The
regularized and constrained budgets are specified with each study. UFGM and
GAGD+ stop at their configured target or the budget limit. The other methods
are evaluated at both reported thresholds along a single run. A target not
reached within the budget is treated as a non-hit.

\subsection{Regression without regularization on synthetic data}
\label{app:setup-details}

\subsubsection{Setup}
\label{app:controlled-instance-details}

The controlled study considers the unregularized objective \(F(\vu)=\norm{\mA\vu-\vb}_p^p/p\) with \(1<p<2\).

\paragraph{Synthetic instance generation.}
The instances have a prescribed spectral norm, condition number, initial distance, and optimal residual scale. Their minimizers and optimal values are known by construction. Algorithm~\ref{alg:controlled-instance-generation} gives the generation procedure.
The matrix \(\mA\in \R^{m\times n}\) has condition number \(\kappa\) and spectral norm \(\norm{\mA}=1\). Let \(\mU \in \R^{m \times n}\) and \(\mV \in \R^{n \times n}\) have orthonormal columns obtained from independent Gaussian matrices, and set
\begin{equation*}
    \mA = \mU \operatorname{Diag}(\sigma_1, \ldots, \sigma_n) \mV^\top,
    \qquad
    \sigma_j = \kappa^{-\frac{j-1}{n-1}},\quad j=1,\ldots,n.
\end{equation*}
The singular values are therefore logarithmically spaced between \(1\) and \(1/\kappa\), which gives \(\norm{\mA}=1\) and \(\kappa(\mA)=\kappa\).

The initial point is \(\vu_0 := \vzero \in \R^n\), and the prescribed minimizer is $\vu^\star = D \vd / \norm{\vd}$, where \(\vd \in \R^n\) is drawn from the standard normal distribution. The initial distance \(D\) is constant across all instances.

The nonzero optimal residual is constructed from a direction in the null space of \(\mA^\top\). For a standard Gaussian vector \(\vg \in \R^m\), set
\begin{equation*}
    \vtheta = \frac{(\mI - \mU\mU^\top)\vg}{\norm{(\mI - \mU\mU^\top)\vg}}.
\end{equation*}
Define \((r_0)_i = \operatorname{sign}(\theta_i) |\theta_i|^{1/(p-1)}\) for each \(i\), and set
\begin{equation*}
    \vr^\star = \frac{\sqrt{m} \rho_{\mathrm{res}}}{\norm{\vr_0}} \vr_0,
    \qquad
    \vb = \mA \vu^\star - \vr^\star,
    \qquad
    \rho_{\mathrm{res}} = 0.2.
\end{equation*}
The vector \(\operatorname{sign}(\vr^\star)\odot|\vr^\star|^{p-1}\) is a positive multiple of \(\vtheta\), so \(\mA^\top\vtheta=0\) gives
\(\nabla F(\vu^\star)=\mA^\top\!\left(\operatorname{sign}(\vr^\star)\odot |\vr^\star|^{p-1}\right)=0\).
Convexity therefore establishes \(\vu^\star\) as a minimizer. The residual root-mean-square magnitude satisfies \(\norm{\vr^\star}/\sqrt m=\rho_{\mathrm{res}}\) for every instance, and the optimal value is
\begin{equation*}
    F_{\text{ref}} = F^\star = \frac{1}{p} \norm{\vr^\star}_p^p.
\end{equation*}

For each seed, the same orthogonal factors \((\mU,\mV)\), optimizer direction, and projected residual direction are reused across all \((p,\kappa)\) pairs. Different seeds use independent random draws. Within each seed and at fixed \(p\), changing \(\kappa\) varies the singular values while holding the orthogonal factors, minimizer, and optimal residual fixed. Changing \(p\) modifies the loss exponent and the corresponding optimal residual.

\begin{algorithm}[t]
\caption{Controlled instance generation for one paired seed}
\label{alg:controlled-instance-generation}
\KwIn{base random seed \(s\)}
\KwInit{
\(m\leftarrow 200\), \(n\leftarrow 50\), \(D\leftarrow 10\),
\(\rho_{\mathrm{res}}\leftarrow 0.2\),
\(\vu_0\leftarrow\vzero\);
\newline
\(p\in\{1.1,1.25,1.5,1.75,1.9\}\),
\(\kappa\in\{1,10,10^2,10^3,10^4\}\)
}

Using seed \(s\), draw independent Gaussian matrices and compute
canonical QR factorizations to obtain
\(\mU\in\R^{m\times n}\) and
\(\mV\in\R^{n\times n}\) with orthonormal columns\;

Draw \(\vd\sim\mathcal N(\vzero,\mI_n)\) and set
\(\vu^\star \leftarrow D\,\vd/\|\vd\|\)\;

Draw \(\vg\sim\mathcal N(\vzero,\mI_m)\) and set
\(\vtheta \leftarrow (\mI-\mU\mU^\top)\vg / \|(\mI-\mU\mU^\top)\vg\|\)\;

\For{\(p\in\{1.1,1.25,1.5,1.75,1.9\}\)}{

    \For{\(i=1,\ldots,m\)}{
        \((r_0)_i \leftarrow \operatorname{sign}(\theta_i)\,|\theta_i|^{1/(p-1)}\)\;
    }

    Rescale the residual to the prescribed RMS magnitude:
    \(\vr^\star \leftarrow \frac{\sqrt m\,\rho_{\mathrm{res}}}{\|\vr_0\|}\,\vr_0\)\;

    Compute the exact reference value \(F_{\mathrm{ref}} \leftarrow \frac{1}{p}\|\vr^\star\|_p^p\)\;

    \For{\(\kappa\in\{1,10,10^2,10^3,10^4\}\)}{

        \For{\(j=1,\ldots,n\)}{
            \(\sigma_j \leftarrow \kappa^{-(j-1)/(n-1)}\)\;
        }

        Construct
        \(\mA \leftarrow \mU\,\operatorname{Diag}(\sigma_1,\ldots,\sigma_n)\,\mV^\top\)\;

        Set \(\vb \leftarrow \mA\vu^\star - \vr^\star\)\;

        Verify the spectral norm, condition number, initial distance,
        residual RMS, stationarity, and primal--dual agreement numerically\;
    }
}
\KwOut{the \(25\) paired controlled instances associated with seed \(s\)}
\end{algorithm}

\paragraph{Curve aggregation and target attainment.}
Convergence curves use a common grid of 200 matrix--vector product budgets: the integers
1 through 20, followed by 180 rounded logarithmically spaced budgets from
21 to 60,000. At each budget, each run contributes its last prescribed output
whose cost does not exceed that budget. Runs that terminate earlier retain
their terminal values at subsequent budgets, so all medians and interquartile
ranges use the ten paired seeds. Quantiles follow the Hyndman--Fan type-7
convention. Aggregation uses these outputs directly, without interpolation
or an additional running minimum. The plotted lines connect the median
values, and markers show the first downward crossings of these lines at
selected gap levels. Target-attainment costs use the first prescribed output
satisfying each target in the full output sequence, independently of the
plotting grid and the displayed vertical range of \(10^{-9}\) to 2.

\subsubsection{Detailed results}
\label{app:results-details}

Table~\ref{tab:controlled-attainment} summarizes target attainment over the
full grid, and Table~\ref{tab:controlled-full-grid} reports the results for
each \((p,\kappa)\) pair. Median first-hit costs are reported together with
attainment counts because they are computed only over runs that reach the
target, and these subsets may differ across methods.

Variable smoothing achieves a lower terminal gap than each of UFGM, AGDA,
and GAGD+ on all \CtrlLowPHighKCount{} paired instances with
\(p\leq1.25\) and \(\kappa\geq100\). At the common product budget, the median
of the per-instance ratio of the baseline terminal gap to the
variable-smoothing terminal gap is \CtrlVvsUFGMRatio{} for UFGM,
\CtrlVvsAGDARatio{} for AGDA, and \CtrlVvsGAGDRatio{} for GAGD+.

\begin{table}[htbp]
    \caption{Target attainment within a budget of $60{,}000$ matrix--vector products. The first panel reports the number of runs reaching each target across all 250 instances. Parentheses give the median first-hit cost among those runs. Each column in the second and third panels reports the number of runs reaching $10^{-9}$ across the 50 instances with a fixed $p$ or $\kappa$. The largest count in each column is bold. Table~\ref{tab:controlled-full-grid} gives the results for individual $(p,\kappa)$ pairs.}
    \label{tab:controlled-attainment}
    \centering
    \setlength{\tabcolsep}{6pt}
    \begin{tabular}{lrr}
        \toprule
        & \multicolumn{2}{c}{Overall hits (median first-hit matvecs)}\\
        \cmidrule(lr){2-3}
        Method & $10^{-6}$ & $10^{-9}$\\
        \midrule
        HS-APG-F & 188 (4,016) & 92 (3,618)\\
        HS-APG-V & \textbf{201 (4,582)} & 107 (9,436)\\
        HS-APG-V-BT & 197 (4,074) & 113 (10,825)\\
        \midrule
        AGDA & 164 (5,258) & 103 (11,494)\\
        UFGM & 193 (3,243) & \textbf{128 (5,064)}\\
        GAGD+ & 187 (3,456) & 108 (6,738)\\
        \bottomrule
    \end{tabular}

    \vspace{0.75em}

    \begin{tabular}{lrrrrr}
        \toprule
        & \multicolumn{5}{c}{$10^{-9}$ hits by exponent $p$}\\
        \cmidrule(lr){2-6}
        Method & $1.1$ & $1.25$ & $1.5$ & $1.75$ & $1.9$\\
        \midrule
        HS-APG-F & 0 & 0 & 18 & 33 & \textbf{41}\\
        HS-APG-V & 0 & 0 & 30 & 36 & \textbf{41}\\
        HS-APG-V-BT & 0 & 0 & \textbf{31} & \textbf{41} & \textbf{41}\\
        \midrule
        AGDA & 0 & 0 & 22 & 40 & \textbf{41}\\
        UFGM & \textbf{20} & \textbf{7} & 28 & 37 & 36\\
        GAGD+ & 5 & 0 & 22 & 40 & \textbf{41}\\
        \bottomrule
    \end{tabular}

    \vspace{0.75em}

    \begin{tabular}{lrrrrr}
        \toprule
        & \multicolumn{5}{c}{$10^{-9}$ hits by $\kappa(\mA)$}\\
        \cmidrule(lr){2-6}
        Method & $1$ & $10$ & $10^2$ & $10^3$ & $10^4$\\
        \midrule
        HS-APG-F & 30 & 28 & 20 & 13 & 1\\
        HS-APG-V & 30 & 30 & \textbf{30} & 16 & 1\\
        HS-APG-V-BT & 30 & 29 & 28 & \textbf{24} & \textbf{2}\\
        \midrule
        AGDA & 30 & 28 & 24 & 20 & 1\\
        UFGM & \textbf{47} & \textbf{40} & 28 & 13 & 0\\
        GAGD+ & 35 & 30 & 22 & 20 & 1\\
        \bottomrule
    \end{tabular}
\end{table}

\FloatBarrier
\begin{figure}[p]
    \centering
    \includegraphics[width=\textwidth]{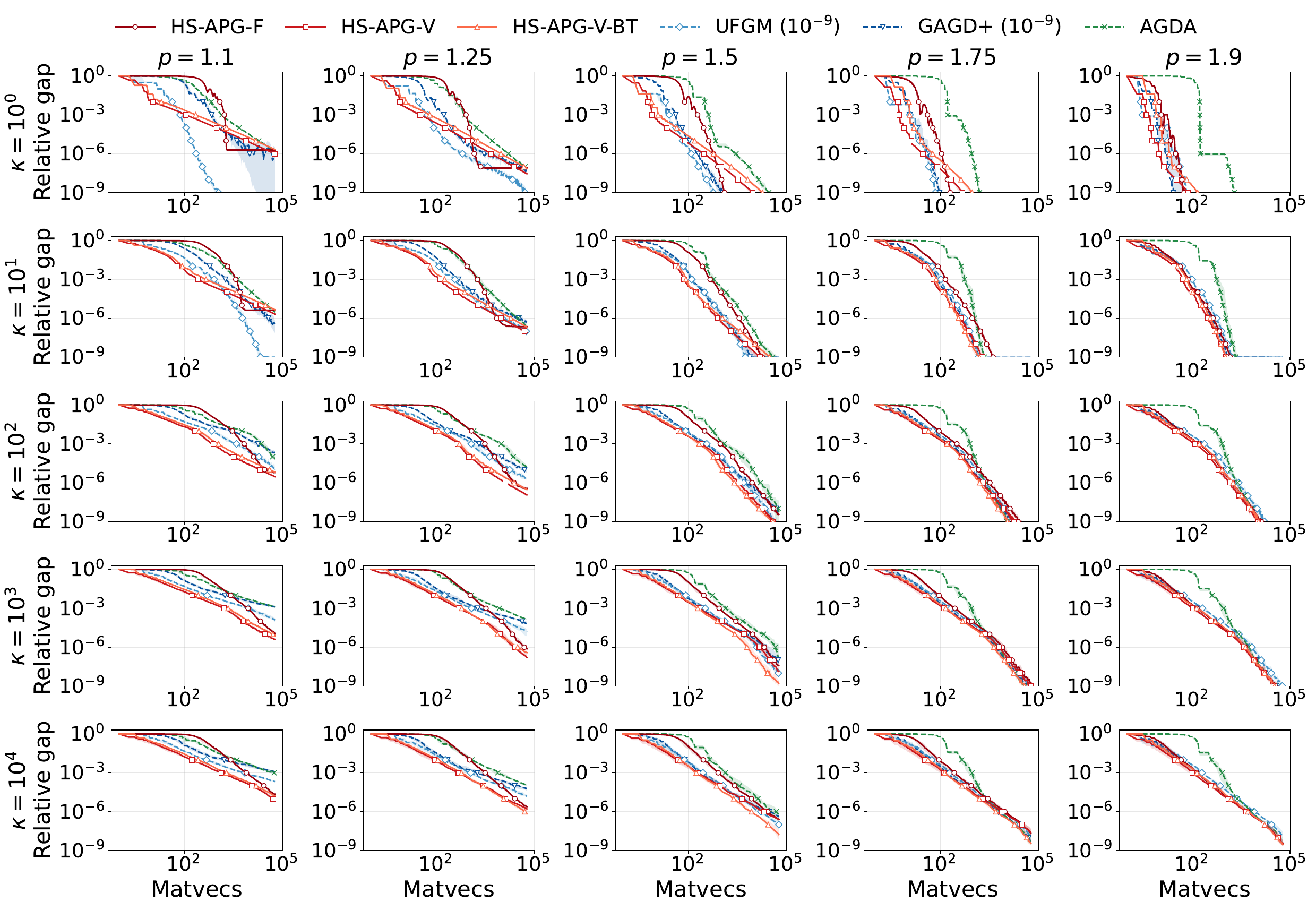}
    \caption{Complete controlled convergence grid over all five exponents
    and five condition numbers. Curves show median relative gaps over ten
    paired seeds against the number of matrix--vector products. Shaded
    bands span the 25th to 75th percentiles. UFGM and GAGD+ use target
    \(10^{-9}\), as in Figure~\ref{fig:controlled-convergence}.}
    \label{fig:controlled-full-grid}
\end{figure}

\FloatBarrier
\begin{table}[htbp]
    \caption{Target attainment for each $(p,\kappa)$ pair within a budget of $60{,}000$ matrix--vector products. Each method entry gives the number of runs reaching the target out of ten paired seeds, followed by the median first-hit cost over those runs. Costs are measured in matrix--vector products; $0$/-- indicates that no run reaches the target. UFGM and GAGD+ are configured for target $10^{-9}$, while the other methods do not require an accuracy target.}
    \label{tab:controlled-full-grid}
    \centering
    \setlength{\tabcolsep}{2pt}
    \begin{tabular}{rrcrrrrrr}
        \toprule
        $p$ & $\kappa$ & Target & F & V & V-BT & AGDA & UFGM & GAGD+\\
        \midrule
        $1.1$ & $1$ & $10^{-6}$ & 0/-- & 8/57,687 & 0/-- & 0/-- & 10/237 & 8/4,753.5\\
        $1.1$ & $1$ & $10^{-9}$ & 0/-- & 0/-- & 0/-- & 0/-- & 10/1,096.5 & 5/12,801\\
        \addlinespace[2pt]
        $1.1$ & $10$ & $10^{-6}$ & 0/-- & 0/-- & 0/-- & 0/-- & 10/5,434.5 & 10/27,205.5\\
        $1.1$ & $10$ & $10^{-9}$ & 0/-- & 0/-- & 0/-- & 0/-- & 10/20,332.5 & 0/--\\
        \addlinespace[2pt]
        $1.1$ & $10^2$ & $10^{-6}$ & 0/-- & 0/-- & 0/-- & 0/-- & 0/-- & 0/--\\
        $1.1$ & $10^2$ & $10^{-9}$ & 0/-- & 0/-- & 0/-- & 0/-- & 0/-- & 0/--\\
        \addlinespace[2pt]
        $1.1$ & $10^3$ & $10^{-6}$ & 0/-- & 0/-- & 0/-- & 0/-- & 0/-- & 0/--\\
        $1.1$ & $10^3$ & $10^{-9}$ & 0/-- & 0/-- & 0/-- & 0/-- & 0/-- & 0/--\\
        \addlinespace[2pt]
        $1.1$ & $10^4$ & $10^{-6}$ & 0/-- & 0/-- & 0/-- & 0/-- & 0/-- & 0/--\\
        $1.1$ & $10^4$ & $10^{-9}$ & 0/-- & 0/-- & 0/-- & 0/-- & 0/-- & 0/--\\
        \addlinespace[2pt]
        \midrule
        $1.25$ & $1$ & $10^{-6}$ & 10/1,391 & 10/4,797 & 10/9,593 & 10/16,126 & 10/550.5 & 10/2,439\\
        $1.25$ & $1$ & $10^{-9}$ & 0/-- & 0/-- & 0/-- & 0/-- & 7/42,237 & 0/--\\
        \addlinespace[2pt]
        $1.25$ & $10$ & $10^{-6}$ & 10/6,508 & 10/9,583 & 10/19,153.5 & 10/32,451.5 & 10/7,606.5 & 10/22,171.5\\
        $1.25$ & $10$ & $10^{-9}$ & 0/-- & 0/-- & 0/-- & 0/-- & 0/-- & 0/--\\
        \addlinespace[2pt]
        $1.25$ & $10^2$ & $10^{-6}$ & 10/22,601 & 10/13,096 & 10/24,017 & 0/-- & 3/49,131 & 0/--\\
        $1.25$ & $10^2$ & $10^{-9}$ & 0/-- & 0/-- & 0/-- & 0/-- & 0/-- & 0/--\\
        \addlinespace[2pt]
        $1.25$ & $10^3$ & $10^{-6}$ & 8/44,925 & 10/25,446 & 10/31,052 & 0/-- & 0/-- & 0/--\\
        $1.25$ & $10^3$ & $10^{-9}$ & 0/-- & 0/-- & 0/-- & 0/-- & 0/-- & 0/--\\
        \addlinespace[2pt]
        $1.25$ & $10^4$ & $10^{-6}$ & 0/-- & 3/44,134 & 7/46,057 & 0/-- & 0/-- & 0/--\\
        $1.25$ & $10^4$ & $10^{-9}$ & 0/-- & 0/-- & 0/-- & 0/-- & 0/-- & 0/--\\
        \addlinespace[2pt]
        \midrule
        $1.5$ & $1$ & $10^{-6}$ & 10/588 & 10/302 & 10/603 & 10/2,908.5 & 10/124.5 & 10/186\\
        $1.5$ & $1$ & $10^{-9}$ & 10/1,256 & 10/9,410 & 10/18,816 & 10/30,009.5 & 10/580.5 & 10/960\\
        \addlinespace[2pt]
        $1.5$ & $10$ & $10^{-6}$ & 10/3,416 & 10/741 & 10/1,110 & 10/5,830.5 & 10/1,129.5 & 10/1,197\\
        $1.5$ & $10$ & $10^{-9}$ & 8/15,521 & 10/15,560 & 9/30,240 & 8/36,559 & 10/5,431.5 & 10/7,426.5\\
        \addlinespace[2pt]
        $1.5$ & $10^2$ & $10^{-6}$ & 10/9,017 & 10/3,604 & 10/2,457 & 10/16,241 & 10/4,914 & 10/6,732\\
        $1.5$ & $10^2$ & $10^{-9}$ & 0/-- & 10/38,734 & 8/37,333.5 & 4/39,737.5 & 8/37,965 & 2/26,748\\
        \addlinespace[2pt]
        $1.5$ & $10^3$ & $10^{-6}$ & 10/21,309 & 10/15,033 & 10/6,272 & 6/40,467 & 10/12,583.5 & 10/20,395.5\\
        $1.5$ & $10^3$ & $10^{-9}$ & 0/-- & 0/-- & 4/53,079 & 0/-- & 0/-- & 0/--\\
        \addlinespace[2pt]
        $1.5$ & $10^4$ & $10^{-6}$ & 10/27,491 & 10/21,173 & 10/8,664 & 8/32,066.5 & 10/16,818 & 9/29,997\\
        $1.5$ & $10^4$ & $10^{-9}$ & 0/-- & 0/-- & 0/-- & 0/-- & 0/-- & 0/--\\
        \addlinespace[2pt]
        \bottomrule
    \end{tabular}
\end{table}

\begin{table}[htbp]
    \centering
    \setlength{\tabcolsep}{2pt}
    Table~\ref{tab:controlled-full-grid} (continued).\par
    \vspace{1ex}
    \begin{tabular}{rrcrrrrrr}
        \toprule
        $p$ & $\kappa$ & Target & F & V & V-BT & AGDA & UFGM & GAGD+\\
        \midrule
        $1.75$ & $1$ & $10^{-6}$ & 10/108 & 10/30 & 10/60 & 10/926 & 10/36 & 10/42\\
        $1.75$ & $1$ & $10^{-9}$ & 10/201 & 10/479 & 10/956.5 & 10/1,596.5 & 10/72 & 10/94.5\\
        \addlinespace[2pt]
        $1.75$ & $10$ & $10^{-6}$ & 10/947 & 10/408 & 10/360 & 10/1,117.5 & 10/550.5 & 10/475.5\\
        $1.75$ & $10$ & $10^{-9}$ & 10/3,467 & 10/1,790 & 10/1,584 & 10/2,269 & 10/1,479 & 10/1,459.5\\
        \addlinespace[2pt]
        $1.75$ & $10^2$ & $10^{-6}$ & 10/3,348 & 10/1,977 & 10/1,630.5 & 10/2,521 & 10/2,497.5 & 10/1,897.5\\
        $1.75$ & $10^2$ & $10^{-9}$ & 10/22,753 & 10/15,310 & 10/10,081 & 10/11,624 & 10/14,859 & 10/13,275\\
        \addlinespace[2pt]
        $1.75$ & $10^3$ & $10^{-6}$ & 10/8,101 & 10/6,454 & 10/4,415.5 & 10/5,752 & 10/6,673.5 & 10/5,172\\
        $1.75$ & $10^3$ & $10^{-9}$ & 3/53,882 & 6/48,726 & 10/34,238 & 10/39,497 & 7/41,529 & 10/40,213.5\\
        \addlinespace[2pt]
        $1.75$ & $10^4$ & $10^{-6}$ & 10/9,523 & 10/7,675 & 10/5,071 & 10/6,629 & 10/7,653 & 10/6,034.5\\
        $1.75$ & $10^4$ & $10^{-9}$ & 0/-- & 0/-- & 1/52,187 & 0/-- & 0/-- & 0/--\\
        \addlinespace[2pt]
        \midrule
        $1.9$ & $1$ & $10^{-6}$ & 10/32 & 10/10 & 10/20 & 10/172 & 10/30 & 10/18\\
        $1.9$ & $1$ & $10^{-9}$ & 10/53 & 10/74 & 10/148 & 10/1,970 & 10/66 & 10/27\\
        \addlinespace[2pt]
        $1.9$ & $10$ & $10^{-6}$ & 10/401 & 10/300 & 10/345 & 10/1,072 & 10/549 & 10/393\\
        $1.9$ & $10$ & $10^{-9}$ & 10/1,384 & 10/1,051 & 10/995 & 10/2,144.5 & 10/1,560 & 10/1,134\\
        \addlinespace[2pt]
        $1.9$ & $10^2$ & $10^{-6}$ & 10/1,737 & 10/1,334 & 10/1,772.5 & 10/2,458.5 & 10/2,662.5 & 10/1,678.5\\
        $1.9$ & $10^2$ & $10^{-9}$ & 10/9,771 & 10/8,891 & 10/9,518.5 & 10/11,119 & 10/15,151.5 & 10/10,203\\
        \addlinespace[2pt]
        $1.9$ & $10^3$ & $10^{-6}$ & 10/4,273 & 10/3,891 & 10/4,364.5 & 10/5,376 & 10/6,879 & 10/4,479\\
        $1.9$ & $10^3$ & $10^{-9}$ & 10/33,957 & 10/30,345 & 10/33,276.5 & 10/35,409.5 & 6/46,504.5 & 10/35,482.5\\
        \addlinespace[2pt]
        $1.9$ & $10^4$ & $10^{-6}$ & 10/4,998 & 10/4,574 & 10/5,058.5 & 10/6,476 & 10/7,944 & 10/5,194.5\\
        $1.9$ & $10^4$ & $10^{-9}$ & 1/52,284 & 1/54,470 & 1/50,542 & 1/58,525 & 0/-- & 1/54,708\\
        \addlinespace[2pt]
        \bottomrule
    \end{tabular}
\end{table}

\FloatBarrier
\subsection{Regression without regularization on real-world data}
\label{app:realdata-experiments}

\subsubsection{Datasets and preprocessing}
\label{app:realdata_descriptions}

We construct regression problems from six UCI datasets: Power, Year, Protein, KEGG, CT-Slices, and Buzz \citep{uciHouseholdPower2006,uciYearPredictionMSD2011,uciProteinCASP2013,uciKeggUndirected2011,uciCTSlices2011,uciBuzzSocialMedia2013}. Table~\ref{tab:realdata-designs} lists the response variables and the resulting matrix dimensions. We exclude Date and Time from Power, the textual Pathway field from KEGG, and \texttt{patientId} from CT-Slices. Protein uses the nine physicochemical variables as predictors. Year uses the first \(463{,}715\) observations in the official ordered training split; the final \(51{,}630\) test observations are excluded. For Buzz, we use \texttt{Twitter.data} from the Twitter regression subset, with the rightmost column as the response.

Rows with missing values in any selected predictor or response are removed before preprocessing statistics are computed. For Power, both \texttt{?} and empty numerical fields are treated as missing. For CT-Slices, the value \(-0.25\), which denotes histogram bins outside the image, is retained. Each predictor and the response are centered and scaled by their population root-mean-square deviations. Constant predictors are removed, and numerical redundancy is assessed by column-pivoted QR without using the response. The tolerance accounts for matrix dimensions and floating-point precision, with a separate tolerance for KEGG to account for the decimal precision of its source data. This procedure removes two redundant predictors from KEGG and five constant and six redundant predictors from CT-Slices. No predictors are removed from the other four datasets. All methods use the same retained predictor sets.

An intercept is then appended, and the augmented design matrix is divided by its spectral norm so that \(\|\mA\|_2=1\). This final scaling is not applied to \(\vb\). No logarithmic transformation or row normalization is used. Table~\ref{tab:realdata-designs} also reports the condition number \(\kappa_2(\mA)=\sigma_{\max}(\mA)/\sigma_{\min}(\mA)\), computed from the processed matrix used by all methods.

\begin{table}[htbp]
    \caption{Processed real-data optimization designs.  The column count
    \(n\) includes the intercept.}
    \label{tab:realdata-designs}
    \centering
    \setlength{\tabcolsep}{6pt}
    \begin{tabular}{llrrr}
        \hline
        Dataset & Response & \(m\) & \(n\) & \(\kappa_2(\mA)\) \\
        \hline
        Power
            & Global active power
            & \(2{,}049{,}280\) & \(7\) & \(3.8417\) \\
        YearPredictionMSD
            & Release year
            & \(463{,}715\) & \(91\) & \(13.260\) \\
        Protein
            & RMSD
            & \(45{,}730\) & \(10\) & \(67.626\) \\
        KEGG
            & Edge count
            & \(64{,}608\) & \(26\) & \(100.90\) \\
        CT-Slices
            & Axial location
            & \(53{,}500\) & \(374\) & \(583.46\) \\
        Buzz/Twitter
            & Mean active discussions
            & \(583{,}250\) & \(78\) & \(7278.3\) \\
        \hline
    \end{tabular}
\end{table}

\subsubsection{High-precision reference construction}
\label{app:realdata_references}

For each dataset and exponent \(p\in(1,2)\), let
\begin{equation*}
\vu^\star\in\arg\min_{\vu}F_p(\vu),
\qquad
F_p^\star=F_p(\vu^\star)
=\min_{\vu}\left\{\frac{1}{p}\norm{\mA\vu-\vb}_p^p\right\}.
\end{equation*}
We compute a primal candidate \(\vu_{\mathrm{ref}}\) using BFGS with Hager--Zhang line search in 512-bit arithmetic \citep{nocedal2006numerical,hzlinesearch2005}. Its objective value gives the upper bound
\begin{equation*}
F_p^{\mathrm U}:=F_p(\vu_{\mathrm{ref}})\geq F_p^\star.
\end{equation*}
We assess reference accuracy through a primal--dual interval \citep{pdrc-2016-icml}. With a lower bound \(F_p^{\mathrm L}\leq F_p^\star\) and the common normalization point \(\vu_0=\vzero\), we require
\begin{equation*}
\frac{F_p^{\mathrm U}-F_p^{\mathrm L}}
     {F_p(\vu_0)-F_p^{\mathrm L}}
\leq 10^{-10}.
\end{equation*}
Objective gaps and stopping targets use \(F_{\mathrm{ref}}=F_p^{\mathrm L}\). The primal candidate \(\vu_{\mathrm{ref}}\) is used only for diagnostics and distance-based parameters.

Let \(q=p/(p-1)\). The Fenchel dual is
\begin{equation*}
\max_{\vlambda:\,\mA^\top\vlambda=0}
\left\{
D(\vlambda):=\vlambda^\top\vb-\frac{1}{q}\norm{\vlambda}_q^q
\right\}.
\end{equation*}
The primal residual gives the direction
\begin{equation*}
\vlambda_0:=-\operatorname{sign}(\vr_{\mathrm{ref}})\odot|\vr_{\mathrm{ref}}|^{p-1},
\qquad
\vr_{\mathrm{ref}}:=\mA\vu_{\mathrm{ref}}-\vb,
\end{equation*}
which satisfies \(\mA^\top\vlambda_0=-\nabla F_p(\vu_{\mathrm{ref}})\). We project this direction onto \(\ker(\mA^\top)\):
\begin{equation*}
\widehat{\vlambda}:=\vlambda_0-\mA(\mA^\top\mA)^{-1}\mA^\top\vlambda_0,
\end{equation*}
with the linear systems solved using a high-precision Gram factorization.

Numerical projection can leave a nonzero feasibility residual \(\mA^\top\widehat{\vlambda}\). To account for this residual, let \(0<\underline{\sigma}(\mA)\leq\sigma_{\min}(\mA)\) be a certified lower bound and set \(R_{\vu}=2\norm{\vb}_p/\underline{\sigma}(\mA)\). Since \(F_p(\vu^\star)\leq F_p(\vzero)\), we have \(\norm{\mA\vu^\star-\vb}_p\leq\norm{\vb}_p\), which implies \(\norm{\vu^\star}_2\leq R_{\vu}\). Fenchel's inequality then gives the corrected lower bound
\begin{equation*}
F_p^{\mathrm L}:=\max_{\alpha\geq0}\left\{
\alpha\left(\widehat{\vlambda}^{\top}\vb-R_{\vu}\norm{\mA^\top\widehat{\vlambda}}_2\right)
-\frac{\alpha^q}{q}\norm{\widehat{\vlambda}}_q^q
\right\}.
\end{equation*}
Thus \(F_p^{\mathrm L}\leq F_p^\star\leq F_p^{\mathrm U}\).

This residual-based construction meets the accuracy requirement for 17 of the 18 problems. For CT-Slices at \(p=1.1\), we instead maximize \(D(Pz)\) using projected L-BFGS with Hager--Zhang line search, where
\begin{equation*}
P:=I-\mA(\mA^\top\mA)^{-1}\mA^\top.
\end{equation*}
The primal candidate remains fixed. The resulting direction is reprojected, and the same corrected lower bound is evaluated in 512-bit arithmetic. This reduces the normalized interval width from approximately \(1.08\times10^{-6}\) to \(2.16\times10^{-11}\). An independent 768-bit computation also satisfies the accuracy requirement, using reversed summation order and the weaker singular-value bound \(17/10000\). All 18 reference intervals satisfy the \(10^{-10}\) requirement.

\subsubsection{Algorithm settings and evaluation}
\label{app:realdata-settings}

\paragraph{Reference-distance setting.}
\label{app:realdata_known_distance}

Fixed smoothing, variable smoothing, and variable smoothing with backtracking start from \(\vu_0=\vzero\) and use \(D_{\mathrm{ref}}=\|\vu_{\mathrm{ref}}\|_2\) as a distance proxy. Fixed smoothing sets \(D=D_{\mathrm{ref}}\), while both variable-smoothing variants set \(R=D_{\mathrm{ref}}\) in their radius formulas. The primal--dual certificate bounds the objective error at \(\vu_{\mathrm{ref}}\); it does not certify \(D_{\mathrm{ref}}\) as an upper bound on the initial or subsequent distances to a minimizer.

\paragraph{Least-squares initialization.}
\label{app:realdata_distance_free}

Variable smoothing with least-squares initialization follows Appendix~\ref{app:ols-variant}. The least-squares point is computed from zero by unpreconditioned LSQR, using only products with \(\mA\) and \(\mA^\top\). LSQR stops when an explicitly computed normal residual satisfies
\[
    \frac{\|\mA^\top(\mA\vu_{\mathrm{LS}}-\vb)\|_2}
         {\|\mA^\top\vb\|_2}\le 10^{-10}.
\]
This residual is checked every ten iterations, when the recurrence indicates convergence or breakdown, and near the budget limit if sufficient products remain.

All products required for initialization count toward the shared budget, including the initial adjoint product and the residual checks. Initialization is performed separately for each problem, with no reuse across exponents. All 18 initializations satisfy the stopping criterion within the budget. The resulting least-squares residuals are nonzero in all 18 problems and determine the smoothing scale used in the subsequent variable-smoothing iterations, without reference-distance information.

Target attainment is evaluated at every completed LSQR iteration and every subsequent accelerated output, using the common denominator \(F(\mathbf 0)-F_{\mathrm{ref}}\). If a target is reached during initialization, its attainment cost is the cumulative cost at that iterate. As for the other methods, products used only to evaluate objectives for reporting are counted separately and excluded from the algorithmic budget. These products are not reused to check the LSQR stopping criterion or to compute the smoothing scale.

\paragraph{Evaluation and cost accounting.}
\label{app:realdata_evaluation}

Variable smoothing with least-squares initialization begins with LSQR at zero; all other methods begin optimization at zero. Each run has a budget of \(60{,}000\) matrix--vector products, including initialization and all line-search products. Product counts and objective evaluations follow Appendix~\ref{app:method-details}.

Target attainment is evaluated at relative gaps \(10^{-4}\) and \(10^{-6}\), using the common denominator \(F(\mathbf 0)-F_{\mathrm{ref}}\). UFGM and GAGD+ are run separately for each target, and each target is evaluated using its matching configuration. The remaining methods are evaluated at both targets along a single run.

The four smoothing variants also have a \(30{,}000\)-iteration limit. Fixed smoothing uses this full horizon to set its radius in \eqref{eq:fixed-radius}, and its first-hit costs at both targets are obtained from the same run.

\subsubsection{Detailed results and adaptive restart}
\label{app:realdata_results}

\paragraph{Target attainment.}
Tables~\ref{tab:realdata-refdist-attainment-1e-4} and~\ref{tab:realdata-refdist-attainment-1e-6} report first-hit costs at the two relative-gap targets.

\begin{table}[htbp]
    \caption{Matrix--vector product counts at first attainment for the 18 real-data problems at relative-gap target $10^{-4}$.}
    \label{tab:realdata-refdist-attainment-1e-4}
    \centering
    \setlength{\tabcolsep}{2pt}
    \renewcommand{\arraystretch}{1.08}
    \begin{tabular*}{\textwidth}{@{\extracolsep{\fill}}llrrrrrrr@{}}
        \toprule
        Dataset & $p$ & HS-APG-F & HS-APG-V & HS-APG-V-BT & GAGD+ & UFGM & AGDA & HS-APG-OLS\\
        \midrule
        Power & 1.1 & 3,440 & \textbf{38} & 108 & 483 & 165 & 785 & 77\\
         & 1.5 & 352 & 36 & 44 & 54 & 129 & 646 & \textbf{9}$^{\dagger}$\\
         & 1.9 & 44 & 30 & 36 & 42 & 69 & 610 & \textbf{9}$^{\dagger}$\\
        \addlinespace[3pt]
        Year & 1.1 & 11,700 & 698 & \textbf{232} & 3,039 & 393 & 1,690 & 315\\
         & 1.5 & 1,914 & 410 & \textbf{176} & 549 & 342 & 1,564 & 241\\
         & 1.9 & 266 & 198 & 239 & 279 & 333 & 1,283 & \textbf{193}\\
        \addlinespace[3pt]
        Protein & 1.1 & 17,762 & 5,332 & \textbf{752} & 4,074 & 1,197 & 6,572 & 2,969\\
         & 1.5 & 3,742 & 1,552 & \textbf{717} & 1,260 & 1,179 & 1,822 & 733\\
         & 1.9 & 692 & 576 & 763 & 738 & 1,017 & 1,475 & \textbf{67}\\
        \addlinespace[3pt]
        KEGG & 1.1 & 8,514 & \textbf{1,406} & 1,989 & -- & 4,017 & 41,801 & 3,789\\
         & 1.5 & 2,498 & 848 & 839 & 936 & 1,257 & 2,473 & \textbf{681}\\
         & 1.9 & 478 & 386 & 451 & 510 & 843 & 1,271 & \textbf{55}$^{\dagger}$\\
        \addlinespace[3pt]
        CT-Slices & 1.1 & 17,988 & 4,990 & \textbf{1,368} & 10,677 & 2,478 & 20,846 & 3,635\\
         & 1.5 & 3,202 & 1,218 & \textbf{813} & 1,524 & 1,254 & 1,620 & 1,159\\
         & 1.9 & 594 & 484 & 675 & 666 & 978 & 969 & \textbf{135}$^{\dagger}$\\
        \addlinespace[3pt]
        Buzz & 1.1 & -- & -- & \textbf{29,412} & -- & -- & -- & --\\
         & 1.5 & 42,554 & 40,058 & 19,332 & 30,012 & 39,084 & -- & \textbf{3,923}\\
         & 1.9 & 14,370 & 14,110 & 17,887 & 18,183 & 24,207 & 52,942 & \textbf{1,085}\\
        \midrule
        \multicolumn{2}{l}{Attained / 18} & 17 & 17 & 18 & 16 & 17 & 16 & 17\\
        \bottomrule
    \end{tabular*}
    \par\smallskip
    \begin{minipage}{\textwidth}
    Bold entries indicate the smallest count among all seven methods in each row, including ties. A dash denotes no attainment within the budget of $60{,}000$ matrix--vector products. Counts include initialization (LSQR for HS-APG-OLS) and all line-search products. Products used solely for objective reporting are excluded. Relative gaps use the common denominator $F(\mathbf{0})-F_{\mathrm{ref}}$. GAGD+ and UFGM are configured for $10^{-4}$. $\dagger$ marks attainment during LSQR initialization.
    \end{minipage}
\end{table}

\begin{table}[htbp]
    \caption{Matrix--vector product counts at first attainment for the 18 real-data problems at relative-gap target $10^{-6}$.}
    \label{tab:realdata-refdist-attainment-1e-6}
    \centering
    \setlength{\tabcolsep}{2pt}
    \renewcommand{\arraystretch}{1.08}
    \begin{tabular*}{\textwidth}{@{\extracolsep{\fill}}llrrrrrrr@{}}
        \toprule
        Dataset & $p$ & HS-APG-F & HS-APG-V & HS-APG-V-BT & GAGD+ & UFGM & AGDA & HS-APG-OLS\\
        \midrule
        Power & 1.1 & 5,900 & \textbf{202} & 363 & 1,509 & 234 & 2,241 & 239\\
         & 1.5 & 724 & \textbf{124} & 148 & 264 & 213 & 2,164 & 143\\
         & 1.9 & 112 & 56 & 71 & 78 & 132 & 918 & \textbf{27}\\
        \addlinespace[3pt]
        Year & 1.1 & 47,044 & 3,814 & \textbf{914} & 9,279 & 1,056 & 6,669 & 2,441\\
         & 1.5 & 7,290 & 1,556 & 536 & 1,734 & 780 & 1,963 & \textbf{511}\\
         & 1.9 & 780 & 538 & 582 & 633 & 774 & 1,646 & \textbf{269}\\
        \addlinespace[3pt]
        Protein & 1.1 & 52,006 & 16,416 & \textbf{2,338} & 27,231 & 2,538 & 16,909 & 5,209\\
         & 1.5 & 7,622 & 3,836 & \textbf{824} & 2,490 & 1,341 & 2,183 & 1,353\\
         & 1.9 & 1,430 & 666 & 875 & 1,524 & 1,167 & 2,516 & \textbf{497}\\
        \addlinespace[3pt]
        KEGG & 1.1 & 36,628 & \textbf{9,418} & 15,896 & -- & -- & -- & 42,499\\
         & 1.5 & 8,858 & 2,896 & \textbf{2,193} & 2,457 & 3,438 & 3,811 & 2,705\\
         & 1.9 & 1,070 & 896 & 1,172 & 1,158 & 1,935 & 1,973 & \textbf{673}\\
        \addlinespace[3pt]
        CT-Slices & 1.1 & -- & 46,724 & \textbf{9,878} & -- & 50,169 & -- & 35,971\\
         & 1.5 & 11,452 & 6,336 & \textbf{2,983} & 7,662 & 4,647 & 23,257 & 3,589\\
         & 1.9 & 2,484 & 2,188 & 2,883 & 2,805 & 4,194 & 6,281 & \textbf{1,151}\\
        \addlinespace[3pt]
        Buzz & 1.1 & -- & -- & -- & -- & -- & -- & --\\
         & 1.5 & -- & -- & -- & -- & -- & -- & --\\
         & 1.9 & 58,312 & -- & -- & -- & -- & -- & \textbf{11,117}\\
        \midrule
        \multicolumn{2}{l}{Attained / 18} & 15 & 15 & 15 & 13 & 14 & 13 & 16\\
        \bottomrule
    \end{tabular*}
    \par\smallskip
    \begin{minipage}{\textwidth}
    Bold entries indicate the smallest count among all seven methods in each row, including ties. A dash denotes no attainment within the budget of $60{,}000$ matrix--vector products. Counts include initialization (LSQR for HS-APG-OLS) and all line-search products. Products used solely for objective reporting are excluded. Relative gaps use the common denominator $F(\mathbf{0})-F_{\mathrm{ref}}$. GAGD+ and UFGM are configured for $10^{-6}$.
    \end{minipage}
\end{table}

\FloatBarrier
\paragraph{Convergence trajectories.}
\label{app:realdata-trajectories}

\begin{figure}[tp]
    \centering
    \includegraphics[width=\textwidth]{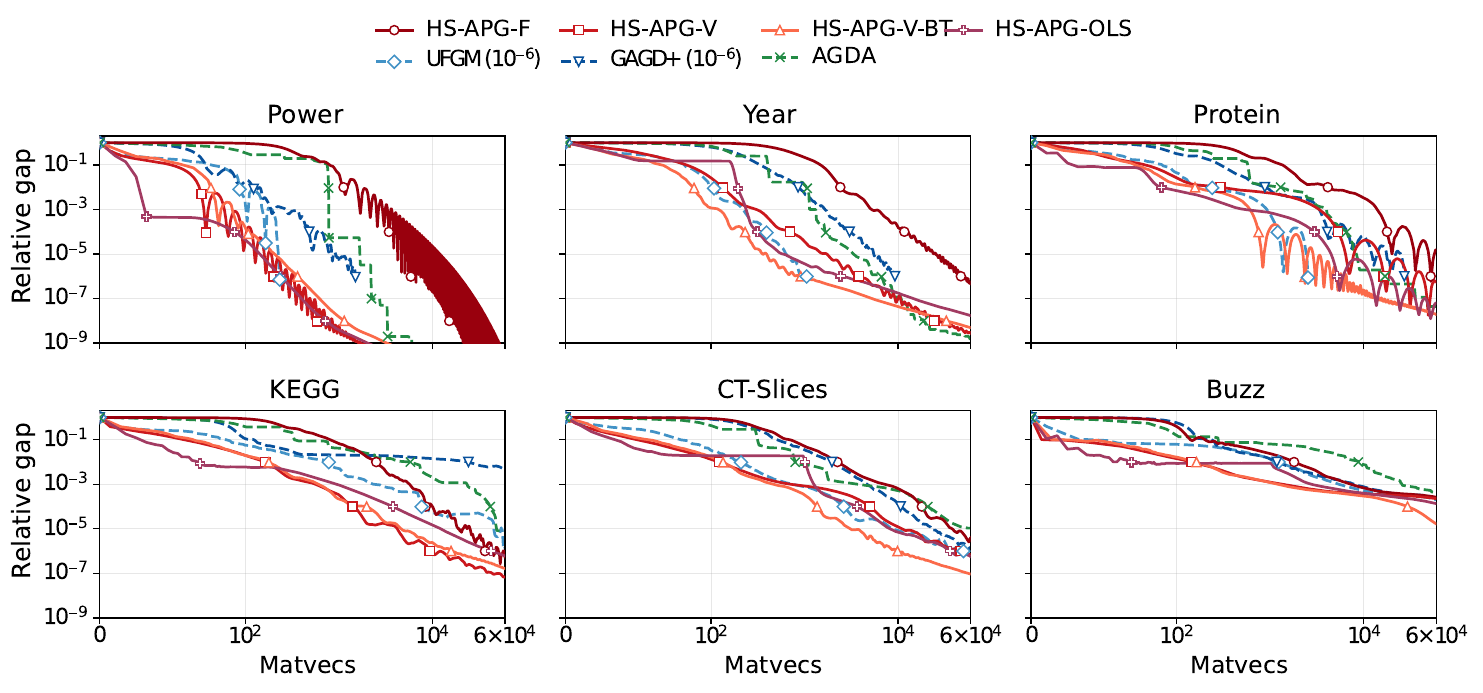}
    \caption{Real-data convergence at \(p=1.1\). Curves show the recorded relative gaps for all seven methods, including LSQR initialization for HS-APG-OLS; UFGM and GAGD+ use target \(10^{-6}\). All panels share axis scales. The matvec axis is linear from 0 to 10 and logarithmic thereafter. Initialization and line-search products are included; products used only for objective reporting are excluded.}
    \label{fig:realdata-convergence-p1p1}
\end{figure}

\begin{figure}[tp]
    \centering
    \includegraphics[width=\textwidth]{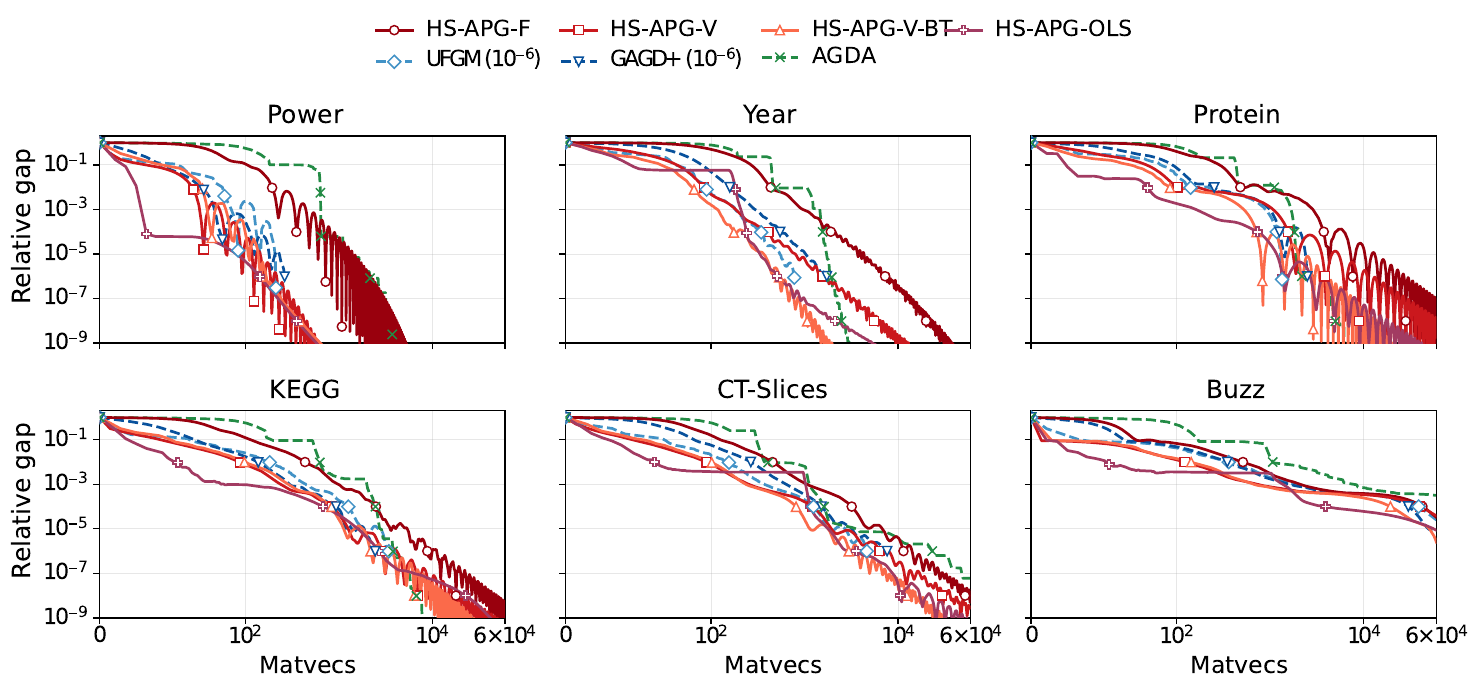}
    \caption{Real-data convergence at \(p=1.5\). Curves show the recorded relative gaps for all seven methods, including LSQR initialization for HS-APG-OLS; UFGM and GAGD+ use target \(10^{-6}\). All panels share axis scales. The matvec axis is linear from 0 to 10 and logarithmic thereafter. Initialization and line-search products are included; products used only for objective reporting are excluded.}
    \label{fig:realdata-convergence}
\end{figure}

\begin{figure}[tp]
    \centering
    \includegraphics[width=\textwidth]{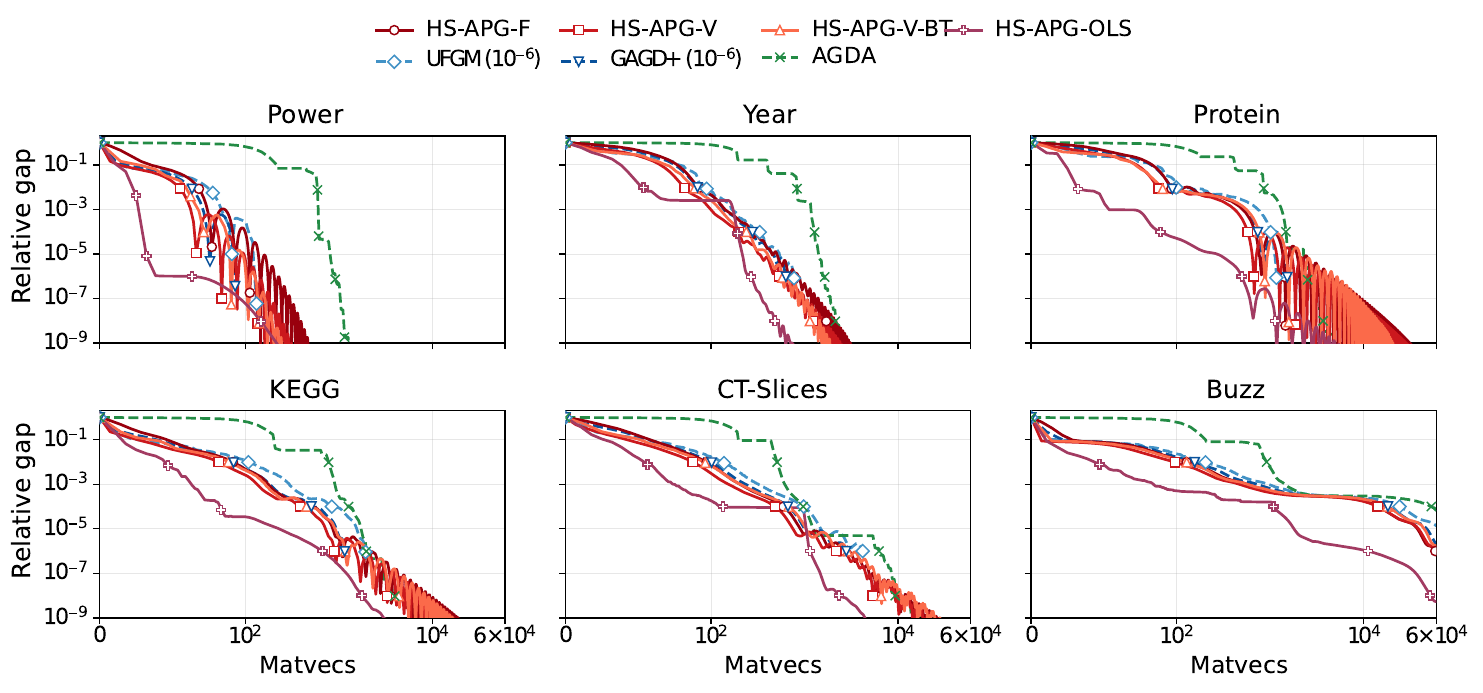}
    \caption{Real-data convergence at \(p=1.9\). Curves show the recorded relative gaps for all seven methods, including LSQR initialization for HS-APG-OLS; UFGM and GAGD+ use target \(10^{-6}\). All panels share axis scales. The matvec axis is linear from 0 to 10 and logarithmic thereafter. Initialization and line-search products are included; products used only for objective reporting are excluded.}
    \label{fig:realdata-convergence-p1p9}
\end{figure}

Figures~\ref{fig:realdata-convergence-p1p1}--\ref{fig:realdata-convergence-p1p9} show convergence on all 18 problems, using every recorded output. For the smoothing methods, the curves report the current accelerated output \(\vy_t\); the least-squares-initialized variant also includes the preceding LSQR iterates. AGDA reports the iterate with the lowest objective value so far, as prescribed by the method. No additional running minimum is applied to any curve. The smoothing trajectories are nonmonotone, so reaching a target once does not imply that subsequent outputs remain below it.

\FloatBarrier
\paragraph{Gradient-restart ablation.}
The convergence curves in Figures~\ref{fig:realdata-convergence-p1p1}--\ref{fig:realdata-convergence-p1p9} exhibit repeated cycles of decrease and rebound. These oscillations motivate adaptive restart, which resets momentum when it appears to oppose descent. We apply the gradient-based rule of \citet{odonoghue2015adaptive} to variable smoothing with and without backtracking.

The rule tests whether \(\langle\vg_t,\vy_t-\vy_{t-1}\rangle>0\), where \(\vg_t=\nabla f_{\rho_t}(\vx_t)\). A positive inner product means that the latest displacement has a component opposite to the negative gradient. The test uses quantities already available in the iteration and requires no additional matrix--vector products, making it inexpensive to add to either method.

The momentum coefficient uses a local iteration counter \(\tau\), with \(\alpha=2/(\tau+1)\), while the smoothing radius follows the global iteration count. After at least five local iterations, a positive test resets \(\vz_t\leftarrow\vy_t\) and \(\tau\leftarrow0\). The current output is retained, as is the accepted curvature estimate in the backtracking variant. Restart therefore does not reject an uphill step or enforce monotonicity.

We compare each method with and without restart on all 18 real-data problems, holding all other settings fixed, including the budget of \(60{,}000\) matrix--vector products. Restart is used only in this ablation. The adaptive rule differs from the horizon-doubling scheme in Appendix~\ref{app:variants} and is not covered by the convergence guarantees established here.

\begin{figure}[tp]
    \centering
    \includegraphics[width=\textwidth]{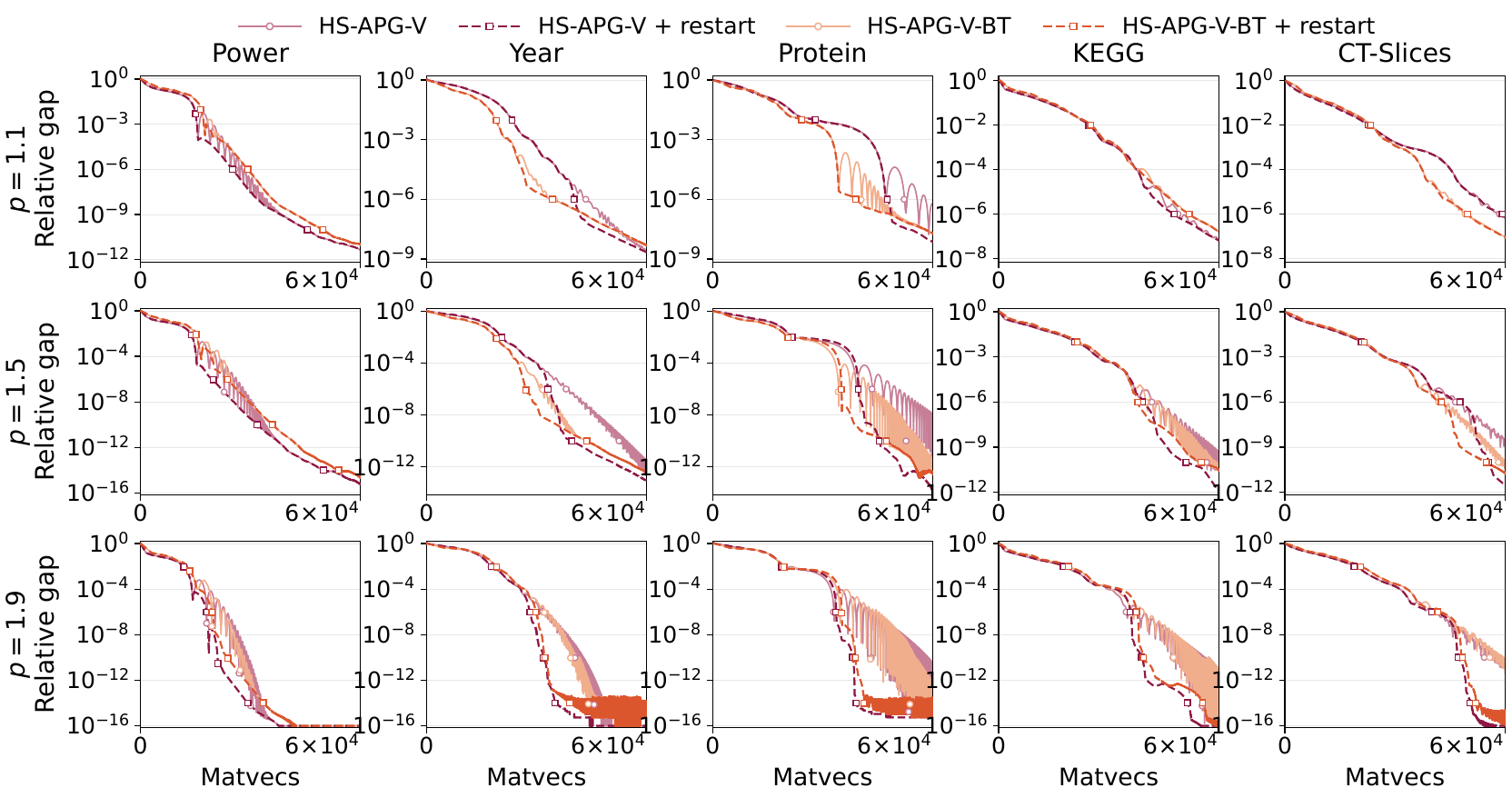}
    \caption{Gradient-restart ablation for HS-APG-V and HS-APG-V-BT. Rows correspond to \(p\in\{1.1,1.5,1.9\}\); columns show the five datasets on which restart is triggered for at least one configuration. Curves retain all recorded outputs, with gaps below \(10^{-16}\) displayed at that level. Vertical limits vary by panel; the matvec axis is linear from 0 to 10 and logarithmic thereafter.}
    \label{fig:realdata-gradient-restart}
\end{figure}

Figure~\ref{fig:realdata-gradient-restart} shows that this simple modification can reduce oscillations and improve terminal accuracy. For variable smoothing on KEGG and CT-Slices at \(p=1.5\), the terminal relative gaps decrease from \(4.79\times10^{-10}\) to \(1.45\times10^{-12}\) and from \(2.88\times10^{-9}\) to \(3.22\times10^{-12}\), respectively. These examples illustrate the value of restart as an inexpensive addition to the smoothing methods when oscillations are pronounced.

Across all 18 problems, restart gives lower, higher, and identical terminal gaps on 10, 4, and 4 problems for variable smoothing, and on 11, 3, and 4 for variable smoothing with backtracking. These counts compare the recorded endpoints directly and include differences near numerical precision, so they do not imply uniform improvement. Buzz remains in these counts but is omitted from the figure: none of its six method--\(p\) configurations triggers restart, and each trajectory is identical to its counterpart without restart.

\FloatBarrier
\subsubsection{Sensitivity to the product budget and smoothing horizon}
\label{app:realdata-budget}

The main experiments use a common budget of \(B_0=60{,}000\) matrix--vector products within the available computational resources. To check whether this cutoff captures the main performance differences, the budget is doubled on four problems selected from the original results. Among the five problems where some methods have not attained \(10^{-6}\), KEGG at \(p=1.1\) has baselines improving near the cutoff, CT-Slices at \(p=1.1\) has both fixed smoothing and GAGD+ close to the target, and Buzz at \(p=1.9\) has several methods approaching it together. Protein at \(p=1.5\), where all seven methods attain the target within \(7{,}622\) products, provides a comparison after early attainment. Selection accounts for additional computational cost and is fixed before the extension; each problem retains all seven methods.

Previously unresolved configurations and all Protein configurations are rerun from their original starting points with a \(120{,}000\)-product cap. The main runs retain the original parameters, including the fixed-smoothing radius set with \(T=30{,}000\), reference bounds, and product-counting rules. Established first-hit costs on the other three problems are reused.

\begin{table}[htbp]
\caption{First-hit costs in the budget-sensitivity study, measured in matrix--vector products. Costs without $\dagger$ or $\ddagger$ are at most 60,000; $\dagger$ denotes $(60{,}000,90{,}000]$ and $\ddagger$ denotes $(90{,}000,120{,}000]$. A dash indicates that the target is not reached within 120,000 products. Superscript $\mathrm{h}$ marks a first-hit cost reused from the original-budget run for that configuration; all other costs are measured in new runs from the original initial points. UFGM and GAGD+ use target-matched configurations. Fixed smoothing retains $T=30{,}000$. On Protein, all seven methods reached $10^{-6}$ within the original budget.}
\label{tab:realdata-budget-hits}
\centering
\setlength{\tabcolsep}{3pt}
\begin{tabular}{lrrrrrrr}
\toprule
Problem & F & V & V-BT & OLS & UFGM & GAGD+ & AGDA\\
\midrule
\multicolumn{8}{c}{Relative-gap target $10^{-4}$}\\
\midrule
KEGG, $p=1.1$ & $8{,}514^{\mathrm{h}}$ & $1{,}406^{\mathrm{h}}$ & $1{,}989^{\mathrm{h}}$ & $3{,}789^{\mathrm{h}}$ & $4{,}017^{\mathrm{h}}$ & -- & $41{,}801$\\
CT-Slices, $p=1.1$ & $17{,}988$ & $4{,}990^{\mathrm{h}}$ & $1{,}368^{\mathrm{h}}$ & $3{,}635^{\mathrm{h}}$ & $2{,}478^{\mathrm{h}}$ & $10{,}677^{\mathrm{h}}$ & $20{,}846$\\
Buzz, $p=1.9$ & $14{,}370^{\mathrm{h}}$ & $14{,}110$ & $17{,}887$ & $1{,}085^{\mathrm{h}}$ & $24{,}207^{\mathrm{h}}$ & $18{,}183^{\mathrm{h}}$ & $52{,}942$\\
\midrule
Protein, $p=1.5$ & $3{,}742$ & $1{,}552$ & $717$ & $733$ & $1{,}179$ & $1{,}260$ & $1{,}822$\\
\midrule
\multicolumn{8}{c}{Relative-gap target $10^{-6}$}\\
\midrule
KEGG, $p=1.1$ & $36{,}628^{\mathrm{h}}$ & $9{,}418^{\mathrm{h}}$ & $15{,}896^{\mathrm{h}}$ & $42{,}499^{\mathrm{h}}$ & $76{,}377^{\dagger}$ & -- & --\\
CT-Slices, $p=1.1$ & $81{,}864^{\dagger}$ & $46{,}724^{\mathrm{h}}$ & $9{,}878^{\mathrm{h}}$ & $35{,}971^{\mathrm{h}}$ & $50{,}169^{\mathrm{h}}$ & $62{,}250^{\dagger}$ & --\\
Buzz, $p=1.9$ & $58{,}312^{\mathrm{h}}$ & $60{,}412^{\dagger}$ & $70{,}042^{\dagger}$ & $11{,}117^{\mathrm{h}}$ & $98{,}910^{\ddagger}$ & $74{,}817^{\dagger}$ & --\\
\midrule
Protein, $p=1.5$ & $7{,}622$ & $3{,}836$ & $824$ & $1{,}353$ & $1{,}341$ & $2{,}490$ & $2{,}183$\\
\bottomrule
\end{tabular}
\end{table}

On the three difficult problems, attainment at \(10^{-4}\) remains 20 of 21 method--problem pairs, while attainment at \(10^{-6}\) increases from 10 to 17 (Table~\ref{tab:realdata-budget-hits}). Four baseline configurations and three smoothing configurations first reach \(10^{-6}\) after \(60{,}000\) products. A smoothing variant has the lowest first-hit cost at both targets on each of the three problems. At \(10^{-6}\), the lowest smoothing costs are \(9{,}418\), \(9{,}878\), and \(11{,}117\) products on KEGG, CT-Slices, and Buzz, compared with \(76{,}377\), \(50{,}169\), and \(74{,}817\) for the fastest baselines. The original budget thus captures the lowest first-hit cost on these problems, while the extension records seven additional hits.

Four additional fixed-smoothing runs set the radius with \(T=60{,}000\), separating a longer run from a change in radius (Table~\ref{tab:realdata-budget-horizon}). This larger horizon delays both target hits on every selected problem and reverses the ordering of fixed smoothing and OLS on KEGG at \(10^{-6}\). The comparison concerns first-hit costs; the theoretical radius is chosen for terminal accuracy.

\begin{table}[htbp]
\caption{Fixed-smoothing first-hit costs when the radius is set with $T=30{,}000$ or $T=60{,}000$, under a maximum budget of $120{,}000$ matrix--vector products. Each value of $T$ produces a separate trajectory, on which both targets are evaluated.}
\label{tab:realdata-budget-horizon}
\centering
\setlength{\tabcolsep}{7pt}
\begin{tabular}{lrrrr}
\toprule
& \multicolumn{2}{c}{Target $10^{-4}$} & \multicolumn{2}{c}{Target $10^{-6}$}\\
\cmidrule(lr){2-3}\cmidrule(lr){4-5}
Problem & $T=30{,}000$ & $T=60{,}000$ & $T=30{,}000$ & $T=60{,}000$\\
\midrule
KEGG, $p=1.1$ & 8,514 & 11,654 & 36,628 & 49,960\\
CT-Slices, $p=1.1$ & 17,988 & 24,594 & 81,864 & 111,764\\
Buzz, $p=1.9$ & 14,370 & 14,876 & 58,312 & 60,370\\
\midrule
Protein, $p=1.5$ & 3,742 & 4,454 & 7,622 & 9,070\\
\bottomrule
\end{tabular}
\end{table}

\begin{figure}[tp]
    \centering
    \includegraphics[width=\textwidth]{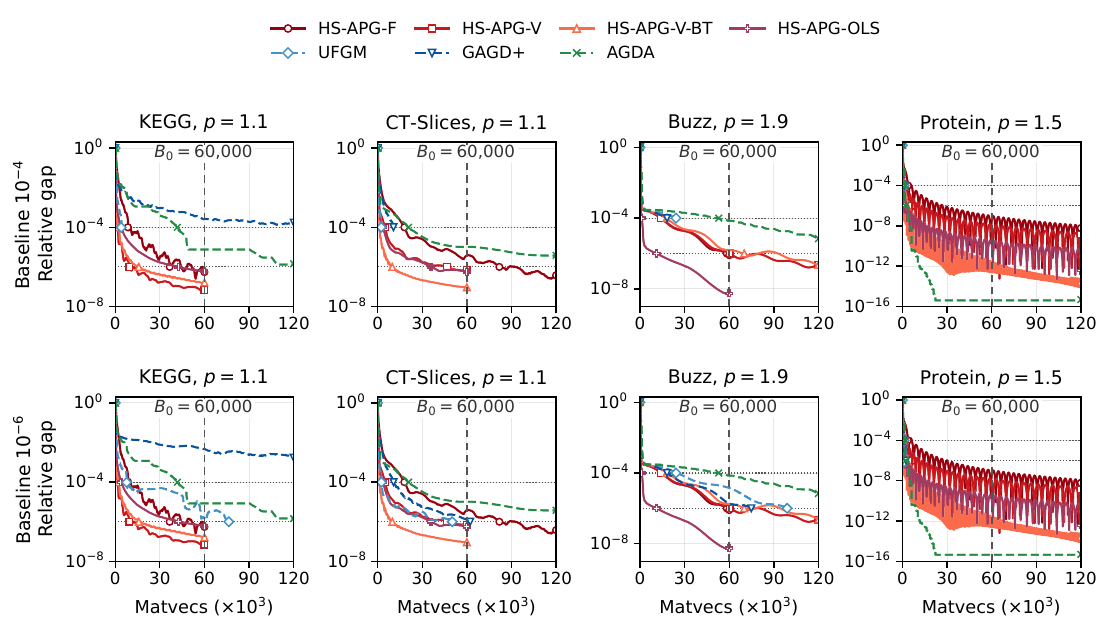}
    \caption{Recorded trajectories under the original parameters. Rows use baseline targets \(10^{-4}\) and \(10^{-6}\); target-independent methods share their trajectory across rows. Vertical dashed lines mark \(B_0=60{,}000\), and horizontal dotted lines mark the targets. The four HS-APG curves on KEGG, V/V-BT/OLS on CT-Slices, and F/OLS on Buzz are historical runs ending at the original budget. Previously unresolved configurations and all Protein configurations use a \(120{,}000\)-product cap. UFGM and GAGD+ stop at the row's target or budget limit. All curves end at their last recorded output; symbols mark initial outputs, first hits, and final outputs. All recorded points are shown, with gaps below \(10^{-16}\) displayed at that level. The x axis is linear in thousands of products; the y axis is logarithmic.}
    \label{fig:realdata-budget-full}
\end{figure}

Protein distinguishes target attainment from subsequent accuracy (Figure~\ref{fig:realdata-budget-full}). V-BT reaches \(10^{-6}\) in 824 products, compared with 2,183 for AGDA. AGDA reaches a lower recorded gap of \(4.91\times10^{-16}\) by 22,207 products, before the original cutoff, and retains it through the extended run. Its curve follows the prescribed best-iterate rule; the smoothing curves show their current outputs.

\par\bigskip

\FloatBarrier
\subsection{Regularized and constrained problems}
\label{app:composite-experiments}

The common algorithm settings and evaluation rules are given in
Appendix~\ref{app:method-details}.

\subsubsection{\texorpdfstring{\(\ell_1\)}{l1}-regularized regression on E2006}
\label{app:regularized_setup}

\paragraph{Setup.}
We use the E2006-tfidf and E2006-log1p training sets, each with \(16{,}087\) observations. Removing all-zero training columns leaves \(150{,}348\) and \(4{,}265{,}669\) features, respectively. Each retained feature is divided by its training \(\ell_2\) norm without centering, while the response is centered and scaled to unit root-mean-square deviation. An intercept is appended, and the augmented design is divided by a conservative spectral-norm upper bound, allowing \(\Lambda=1\) in the algorithmic parameters. All preprocessing statistics are computed from the training set.

We solve~\eqref{eq:lp-regression} with \(p=1.5\) and \(\omega(\vu)=\lambda\norm{\vw}_1\), where \(\vu=(\vw,c)\), \(\vw\) contains the feature coefficients, and \(c\) is the unpenalized intercept. After normalization, the augmented design has the form \(\mA=[\mA_w,s\vone]\), where \(s\) is the intercept-column scale and the prediction intercept is \(sc\). Set \(c_0=\argmin_c f(\mathbf0,c)\) and \(\vr_0=sc_0\vone-\vb\). The smallest penalty for which \((\mathbf0,c_0)\) is optimal is
\begin{equation}
\label{eq:e2006-lambda-max}
\lambda_{\max}=\norm{\mA_w^\top\nabla\psi(\vr_0)}_\infty,
\qquad \nabla\psi(r)=\operatorname{sign}(r)|r|^{1/2},
\end{equation}
where \(\nabla\psi\) acts coordinatewise. We use \(\lambda/\lambda_{\max}\in\{0.01,0.03,0.1,0.3\}\) and initialize every method at \(\vu_0=(\mathbf0,c_0)\). The eight problems and six configurations give \(48\) runs.

\paragraph{Algorithm settings.}
Each run has a budget of \(30{,}000\) matrix--vector products for TFIDF or \(12{,}000\) for log1p. Variable smoothing, with and without backtracking, runs until the budget is exhausted. UFGM and GAGD+ stop at their configured target or budget, whichever comes first. The baseline parameters and common evaluation rules are given in Appendix~\ref{app:method-details}.

\paragraph{Reference construction.}
For reference computation, we combine identical feature columns and sum their coefficients. This preserves predictions and cannot increase the \(\ell_1\) penalty. Conversely, assigning each reduced coefficient to one original column preserves the objective value, so the reduced and original problems have the same optimal value. The compared methods use the original design matrices. We use dual augmented-Lagrangian refinement for lower-penalty problems and high-precision refinement for candidates with small supports. For \(p=1.5\), the Fenchel dual is
\begin{equation*}
\max_{\vv}\;
-\vb^\top\vv-\frac13\norm{\vv}_3^3,
\qquad
\norm{\mA_w^\top\vv}_\infty\leq\lambda,
\quad \vone^\top\vv=0,
\end{equation*}
We compute primal upper bounds \(U\) and dual lower bounds \(L\) on the original problems using 256-bit directed rounding. The lower bound includes a correction for the residual in the intercept equality constraint. Every reference satisfies \((U-L)/(F(\vu_0)-L)\leq10^{-10}\); the largest normalized width is less than \(6.56\times10^{-13}\). We use \(F_{\mathrm{ref}}=L\). Table~\ref{tab:e2006-grid-references} reports the reference values and interval widths.

\paragraph{Distance proxy.}
As in the unregularized experiments, we use the distance from the initial point to a reference candidate. Since E2006 contains many identical feature columns, we first preprocess the reference candidate by distributing each group's total coefficient equally across its columns. Let \(\vu_{\mathrm{ref}}=(\vw_{\mathrm{ref}},c_{\mathrm{ref}})\) be the reference candidate before redistribution. For each group \(G_g\) of identical feature columns, define \(k_g=|G_g|\) and \(z_g=\sum_{j\in G_g}(w_{\mathrm{ref}})_j\). The resulting distance from \(\vu_0\) is
\begin{equation*}
D_{\mathrm{proxy}}=\left(\sum_g\frac{z_g^2}{k_g}+(c_{\mathrm{ref}}-c_0)^2\right)^{1/2}.
\end{equation*}
We use this distance as an empirical estimate, setting \(R=D_{\mathrm{proxy}}\) for variable smoothing with and without backtracking.

\begin{table}[htbp]
\centering
\caption{E2006 reference lower bounds $F_{\mathrm{ref}}$, normalized interval widths, and distance proxies $D_{\mathrm{proxy}}$. All eight normalized widths are below $10^{-10}$.}
\label{tab:e2006-grid-references}
\begin{tabular}{lrrrr}
\toprule
Representation & $\lambda/\lambda_{\max}$ & $F_{\mathrm{ref}}$ & Certified width & $D_{\mathrm{proxy}}$ \\
\midrule
TFIDF & 0.01 & 1683.142621 & $4.7547\times10^{-13}$ & 39293.7320 \\
TFIDF & 0.03 & 3477.862255 & $6.5554\times10^{-13}$ & 50526.3996 \\
TFIDF & 0.1 & 5096.030672 & $2.2936\times10^{-16}$ & 51425.7842 \\
TFIDF & 0.3 & 7567.517556 & $6.0878\times10^{-16}$ & 24027.6634 \\
log1p & 0.01 & 1334.663936 & $4.8745\times10^{-13}$ & 47102.5470 \\
log1p & 0.03 & 3239.009632 & $4.2870\times10^{-13}$ & 54585.3624 \\
log1p & 0.1 & 5159.912785 & $2.3311\times10^{-16}$ & 54754.1167 \\
log1p & 0.3 & 7525.174126 & $1.1840\times10^{-15}$ & 10829.1211 \\
\bottomrule
\end{tabular}
\end{table}

\paragraph{Results.}
Table~\ref{tab:e2006-grid-terminal} reports the terminal relative gaps. All \EThreeCompleted{} runs exhaust their matrix--vector budgets; none reaches a relative gap of \(10^{-4}\) or \(10^{-6}\). Variable smoothing attains the smallest terminal gap on five problems: TFIDF at ratios \(0.01\), \(0.03\), and \(0.1\), and log1p at \(0.01\) and \(0.03\). Variable smoothing with backtracking attains the smallest gap on the remaining three. Across all configurations, terminal relative gaps range from \(4.6715\times10^{-4}\) to \(8.4810\times10^{-2}\).

We measure the sparsity of the terminal feature coefficients by
\begin{equation*}
s_\eta(\vw)=\#\left\{j:|w_j|>\eta\max\{1,\norm{\vw}_\infty\}\right\}.
\end{equation*}
Table~\ref{tab:e2006-grid-sparsity} reports \(s_{10^{-8}}\), with \(s_{10^{-6}}\) as a sensitivity check. For variable smoothing, the TFIDF counts decrease from \(62{,}550\) to \(16{,}668\), \(830\), and \(267\) as the penalty ratio increases. The corresponding log1p counts are \(343{,}331\), \(69{,}985\), \(2{,}255\), and \(1{,}635\).

\begin{table}[htbp]
\centering
\setlength{\tabcolsep}{6pt}
\caption{Thresholded active feature counts at termination for all 48 runs. Columns correspond to penalty ratios $\lambda/\lambda_{\max}$. Entries report $s_{10^{-8}}$, with $s_{10^{-6}}$ in parentheses, where $s_\eta=\#\{j:|w_j|>\eta\max(1,\|\vw\|_\infty)\}$.}
\label{tab:e2006-grid-sparsity}
\begin{tabular}{lrrrr}
\toprule
\multicolumn{5}{c}{TFIDF} \\
Method / target & $0.01$ & $0.03$ & $0.1$ & $0.3$ \\
\midrule
HS-APG-V & 62,550 (32,069) & 16,668 (9,149) & 830 (349) & 267 (105) \\
HS-APG-V-BT & 67,873 (34,953) & 16,945 (9,779) & 823 (382) & 242 (123) \\
UFGM ($10^{-4}$) & 135,355 (72,783) & 44,494 (12,931) & 2,834 (625) & 609 (192) \\
UFGM ($10^{-6}$) & 135,324 (72,702) & 44,459 (12,926) & 2,825 (624) & 609 (192) \\
GAGD+ ($10^{-4}$) & 118,513 (47,599) & 34,533 (9,924) & 2,209 (424) & 544 (146) \\
GAGD+ ($10^{-6}$) & 118,513 (47,599) & 34,533 (9,924) & 2,209 (424) & 544 (146) \\
\midrule
\multicolumn{5}{c}{log1p} \\
Method / target & $0.01$ & $0.03$ & $0.1$ & $0.3$ \\
\midrule
HS-APG-V & 343,331 (97,795) & 69,985 (23,207) & 2,255 (726) & 1,635 (606) \\
HS-APG-V-BT & 329,372 (109,121) & 62,710 (24,906) & 1,870 (711) & 1,212 (510) \\
UFGM ($10^{-4}$) & 3,105,692 (1,020,374) & 564,126 (97,112) & 14,432 (2,241) & 3,840 (934) \\
UFGM ($10^{-6}$) & 3,107,489 (1,021,937) & 562,944 (96,912) & 14,471 (2,246) & 3,838 (934) \\
GAGD+ ($10^{-4}$) & 2,310,619 (527,518) & 325,926 (48,510) & 8,217 (1,210) & 3,161 (621) \\
GAGD+ ($10^{-6}$) & 2,312,722 (528,050) & 325,926 (48,510) & 8,217 (1,210) & 3,161 (621) \\
\bottomrule
\end{tabular}
\end{table}

\FloatBarrier
\subsubsection{Box-constrained image deblurring}
\label{app:constrained_setup}

\paragraph{Setup.}
We use the \(512\times512\) \emph{camera} image from scikit-image v0.25.2 \citep{scikit-image}. Averaging nonoverlapping \(4\times4\) blocks and scaling the intensities to \([0,1]\) gives a \(128\times128\) image \(\vu_{\mathrm{true}}\). The observation is \(\vb=\mA\vu_{\mathrm{true}}+\vxi\), where \(\mA\) is periodic convolution with a normalized \(5\times5\) Gaussian kernel of standard deviation \(0.6\) pixels. The noise entries are independent Gaussian variables with standard deviation \(0.01\), generated with seed \(20260919\).

We solve~\eqref{eq:lp-regression} with \(p=1.5\), \(m=n=128^2\), and \(\omega=\iota_{[0,1]^n}\), whose proximal operator is coordinatewise projection onto the box. All methods start from \(\vu_0=\frac12\vone\). Fixed smoothing uses the geometric initial-distance bound \(D=\sqrt n/2=64\), while variable smoothing with and without backtracking use the box diameter \(R=\sqrt n=128\). The spectral-norm bound is \(\Lambda=1\). Each run has a budget of \(20{,}000\) matrix--vector products, and fixed smoothing uses the nominal horizon \(T=10{,}000\). The smoothing methods stop at relative gap \(10^{-6}\), while UFGM and GAGD+ stop at their respective configured targets.

\paragraph{Reference construction.}
Projected Newton continuation produces a feasible reference candidate. For any \(\vv\in\R^m\), the box-constrained problem admits the lower bound
\begin{equation*}
L(\vv)=-\frac13\norm{\vv}_3^3-\vb^\top\vv
       +\sum_{j=1}^n\min\{0,(\mA^\top\vv)_j\}.
\end{equation*}
The final term is the minimum of \(\langle\mA^\top\vv,\vu\rangle\) over \([0,1]^n\). Evaluating the candidate objective and the dual bound using 256-bit directed rounding gives an interval containing the optimal objective value. Its normalized width is \(2.11\times10^{-11}\), below the common \(10^{-10}\) threshold. We use the lower endpoint as \(F_{\mathrm{ref}}\).

\begin{figure}[tp]
    \centering
    \includegraphics[width=\textwidth]{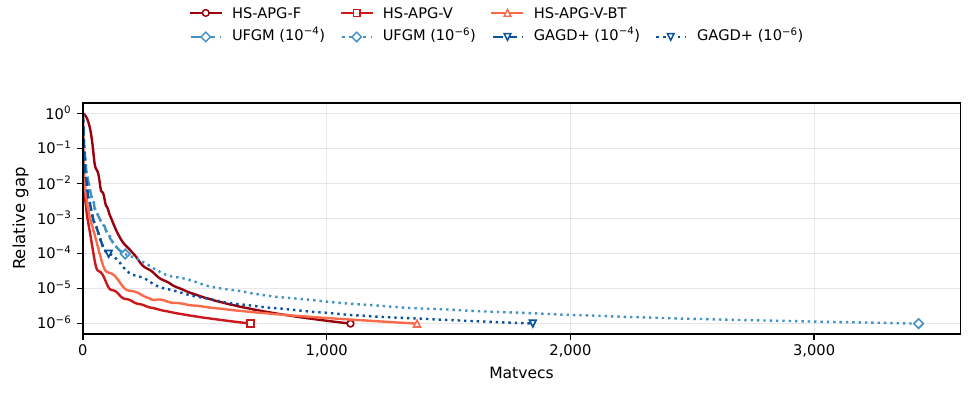}
    \caption{Convergence for box-constrained image deblurring. Curves show relative gaps against matrix--vector products, including line-search trials. Each curve ends at its stopping target. UFGM and GAGD+ configurations are labeled by target.}
    \label{fig:box-convergence}
\end{figure}

\begin{figure}[tp]
    \centering
    \includegraphics[width=\textwidth]{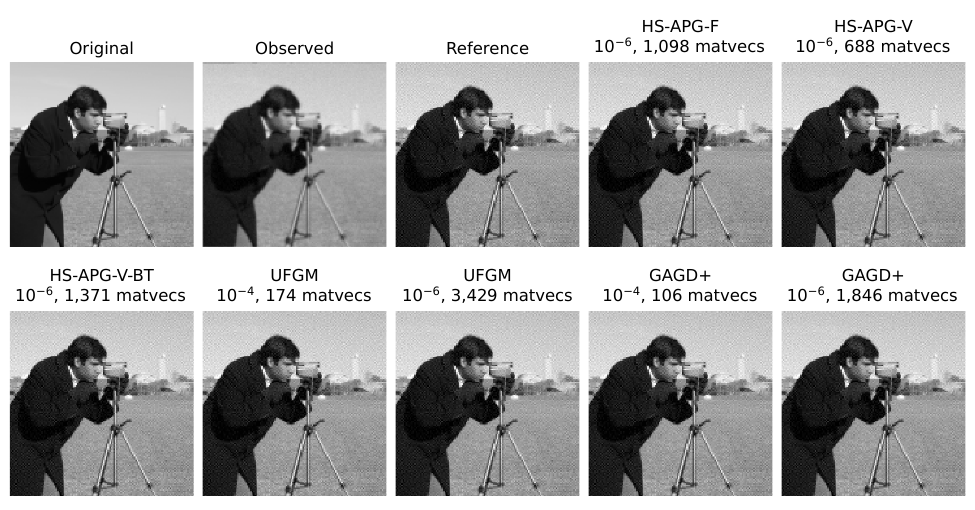}
    \caption{Original image, simulated observation, reference candidate, and terminal reconstructions for box-constrained image deblurring. Method labels give the stopping target and matrix--vector product count. All panels use the grayscale range \([0,1]\).}
    \label{fig:box-reconstructions}
\end{figure}

\paragraph{Results.}
Figure~\ref{fig:box-convergence} shows the convergence curves, and Table~\ref{tab:box-attainment} reports the cost of reaching each target. Variable smoothing uses the fewest matrix--vector products at both targets, reaching \(10^{-4}\) and \(10^{-6}\) in \(44\) and \(688\) products, respectively, compared with \(174\) and \(3{,}429\) for UFGM and \(106\) and \(1{,}846\) for GAGD+. All seven configurations reach their respective stopping targets. Figure~\ref{fig:box-reconstructions} displays the original image, simulated observation, reference candidate, and terminal reconstructions.

\begin{table}[htbp]
\caption{Matrix--vector products required to reach each target in box-constrained image deblurring. UFGM and GAGD+ use the configuration matching each target. Bold entries indicate the smallest count in each column.}
\label{tab:box-attainment}
\centering
\setlength{\tabcolsep}{12pt}
\begin{tabular}{lrr}
\toprule
Method & Target $10^{-4}$ & Target $10^{-6}$ \\
\midrule
HS-APG-F & 208 & 1,098 \\
HS-APG-V & \textbf{44} & \textbf{688} \\
HS-APG-V-BT & 72 & 1,371 \\
UFGM & 174 & 3,429 \\
GAGD+ & 106 & 1,846 \\
\bottomrule
\end{tabular}
\end{table}

\end{document}